\documentclass{article}
\usepackage{graphicx}  
\usepackage{amsmath}
\usepackage{amssymb}
\usepackage{amsthm}
\usepackage{mathtools}
\usepackage{mathrsfs} 
\usepackage[usenames,dvipsnames]{xcolor}
\usepackage{bbm, dsfont}
\usepackage{todonotes}
\usepackage[a4paper]{geometry}
\usepackage[colorlinks=true,linkcolor=NavyBlue,urlcolor=RoyalBlue,citecolor=PineGreen,bookmarks=true,bookmarksopen=true,bookmarksopenlevel=2,unicode=true,linktocpage,
           pdftitle={Loop soup representation of zeta-regularised determinants and covariant Symanzik identities},pdfauthor={P. Perruchaud and I. Sauzedde}]{hyperref}
\usepackage[normalem]{ulem}

\title{Loop soup representation of $\zeta$-regularised determinants\\and covariant Symanzik identities}
\date{February 2024}

\DeclareMathOperator{\Tr}{Tr}
\DeclareMathOperator{\tr}{tr}
\DeclareMathOperator{\Det}{det}
\DeclareMathOperator{\Ker}{ker}
\DeclareMathOperator{\Id}{Id}
\DeclareMathOperator{\End}{End}

\DeclareMathOperator{\Vol}{Vol}
\DeclareMathOperator{\vol}{vol}

\DeclareMathOperator{\grad}{grad}
\DeclareMathOperator{\Range}{Range}

\def\d{\,{\rm d}}

\newtheorem{theorem}{Theorem}[section]

\newtheorem{lemma}[theorem]{Lemma}

\newtheorem{prop}[theorem]{Proposition}
\newtheorem{corollary}[theorem]{Corollary}

\newtheorem{definition}[theorem]{Definition}

\theoremstyle{definition}

\newtheorem{remark}[theorem]{Remark}

\begin{document}

\author{%
P. Perruchaud%
\footnote{Université du Luxembourg, Unité de Recherche en Mathématiques, Maison du Nombre, 6 avenue de la Fonte, L-4364 Esch-sur-Alzette, Grand Duché du Luxembourg
---
Email: \href{mailto:Pierre Perruchaud <pierre.perruchaud@uni.lu>}{pierre.perruchaud@uni.lu}},
I. Sauzedde%
\footnote{University of Warwick, Department of Statistics, Coventry, CV4 7AL, United Kingdom
---
Email: \href{mailto:Isao Sauzedde <isao.sauzedde@warwick.ac.uk>}{isao.sauzedde@warwick.ac.uk}}
}

\maketitle

\begin{abstract}

We derive a stochastic representation for determinants of Laplace-type operators on vectors bundles over manifolds. Namely, inverse powers of those determinants are written as the expectation of a product of holonomies defined over Brownian loop soups. Our results hold over compact manifolds of dimension $d\in\{2,3\}$, in the presence of a mass or a boundary. We derive a few consequences, including some continuity of these determinants as a function of the operator, and the conformal invariance of the determinant on surfaces.

This expression allows us to construct the scalar field minimally coupled to a prescribed random smooth gauge field, which we prove obeys the so-called Symanzik identities. Some of these results are continuous analogues of the work of A.~Kassel and T. Lévy in the discrete.
\smallskip

\smallskip

\end{abstract}
\tableofcontents

\section{Introduction}

\paragraph*{Determinants of Laplace operators}
Determinants of Laplace operators are central objects in many areas of mathematics and physics. In differential topology, they are the building blocks of the so-called analytic torsion \cite{Cheeger,Muller,RaySinger2,RaySinger}, a powerful invariant of manifolds that can discriminate between non-homeomorphic manifolds of the same homotopy type, and is particularly useful in odd dimension. They arise as the partition function of models in various subfields of physics, for instance in gauge theories \cite{Sheffield2,Frohlich,KasselLevy,Symanzik}, or in string theory and Liouville quantum gravity \cite{Sheffield,DHokerPhong,Polyakov}. In statistical mechanics, they are used as weights for probabilities of random maps, which is akin to weights of combinatorial nature related to trees and loops on these maps, or to weights related to discrete Gaussian free fields \cite{Sheffield,Dubedat,KasselKenyon,Kenyon}.

The very definition of those determinants may come with a non-trivial amount of machinery, and a certain part of interpretation. As opposed to discrete, finite models, where Laplace operators are genuine linear operators between finite dimensional spaces, the combinatorial definition is bound to fail in the infinite-dimensional setting, and even the product of the eigenvalues is extremely divergent, since they are far from being perturbations of the identity. On manifolds, which is the case we have in mind, the most common procedure to define them is the so-called $\zeta$-regularisation method, which can be summarised as follows. Let $(\lambda_i)_{i\geq0}$ be the collection of the (positive, with our sign convention) eigenvalues of a Laplace operator $L$, listed with multiplicity. Then formally, one should have on the one hand
\[
\log(\Det(L))= \log\Big(\prod_i \lambda_i\Big)= \sum_i \log(\lambda_i),
\]
and on the other hand, setting $\zeta(s)\coloneqq \sum_i \lambda_i^{-s}$,
\[ 
-\zeta'(0)= -\sum_i \frac{\d}{\d s}_{|s=0} \,e^{-s\log(\lambda_i)}
 =\sum_i \log(\lambda_i),
\]
which suggests defining $\Det(L)$ as $\exp(-\zeta'(0))$. It is important to note that the intermediate steps contain sums and products which are either infinite or ill-defined, so the above is merely motivational. However, it is a striking fact that the function $\zeta$ is well-defined and analytic over a right half plane corresponding to complex numbers with large real part, and that it is the restriction of a (necessarily unique) meromorphic function over $\mathbb C$ that is finite around zero. It is now natural to call this extension $\zeta$, and the determinant is defined by making the leap of faith of setting $\Det(L)=\exp(-\zeta'(0))$.
\medskip

\paragraph*{Alternative representation through the Brownian loop soup}
The study of the determinant often proves to be difficult in the formulation given so far. Indeed, the information extracted from the spectrum must be precise enough to still yield meaningful data after the processes of analytical continuation and removal of singularities.

As a first main result, we derive the new expression \eqref{eq:detintro} below for these determinants, as the expectation of some infinite product taken over the loops of a Brownian loop soup. The formula holds for a large variety of Laplace-type operators in dimension $d\leq 3$.
Some of the properties of the determinant can be deduced almost trivially from this stochastic representation (see the list on page 4), while the initial definition by $\zeta$-regularisation makes those computations much move involved.

We now introduce our framework. We consider a compact connected Riemannian manifold $(M,g)$ with dimension $d\in \{2,3\}$, with or without boundary $\partial M$. Above $M$, we fix a vector bundle $E$, either real or complex, together with a bundle metric $h_E$. Connections on $E$ will be assumed to be compatible with that metric. A function $m\in \mathcal{C}^\infty(M,\mathbb{R}_+)$, often called a mass term, is also fixed, assumed not to be identically vanishing in the case $\partial M=\emptyset$ (this is to ensure that there is no zero eigenvalue for the massive Dirichlet Laplacian on $M$). Then, given a connection $\nabla$ on $E$, one can form a non-negative Laplace operator $L_\nabla = -\frac{1}{2}\Tr(\nabla^2)+m$. It has a well-defined $\zeta$-regularised determinant $\Det(L_\nabla)$, described in a previous paragraph. On the other hand, we can build on $M$ a Brownian loop soup $\mathcal{L}$ with mass $m$, with extinction at the boundary, and intensity $\operatorname{rk}(E)$. For each loop $\ell\in \mathcal{L}$, we can form a scalar $\tr \mathcal Hol^\nabla (\ell) $, where the holonomy $\mathcal Hol^\nabla (\ell)$ is an orthogonal or a unitary endomorphism of the fibre $E_{\ell(0)}$, uniquely determined by $\ell$ and $\nabla$, and the trace is normalised so that the trace of the identity is equal to $1$. The product of these unitary complex numbers, over all the infinitely many loops in the Brownian loop soup, is not an absolutely convergent one. Yet and as we will see, it is possible to define the infinite product, roughly speaking as the limit of an $L^2$-bounded martingale. 

Our formula describes the variation of the determinant as a function of the connection in the precise form
\begin{equation}
\label{eq:detintro}
\det(L_\nabla)^{-1}=C^{\operatorname{rk}E}\,\mathbb{E}\bigg[\prod_{\ell\in \mathcal{L}}  \tr \mathcal Hol^\nabla (\ell)\bigg].
\end{equation}
The constant $C$ depends upon $(M,g)$ but not upon the bundle data $(E, h_E, \nabla)$. For practical purposes, it is not the determinant of $L_\nabla$ that matters, but rather the ratio of two such determinants, so this constant is irrelevant when the metric is fixed.%
\footnote{As far as we understand them, such determinants should be thought of as being defined only up to a multiplicative constant whose choice is a question of convention. Although we do not have hard evidence in this direction, let us mention the works \cite{Dubedat,KasselLevy,AHKHK} where the central objects are indeed such ratios.}
Furthermore, $C$ itself is the inverse determinant of the scalar Laplacian, and the way it depends upon the metric is already well understood by the so-called Polyakov--Alvarez formula, see e.g. \cite[Equation (1.17)]{Osgood}.

The restriction to $d<4$ is not superficial: unless $\nabla$ is a flat connection, the product $\prod_{\ell\in \mathcal{L}} \tr \mathcal Hol^\nabla(\ell)$ can only by defined for $d<4$, since the aforementioned martingale is otherwise overwhelmed by the small loops contribution and not convergent anymore.
Even on average, in dimension $d>4$, this product should as $\exp(\delta^{2-\frac{d}{2}})$, where $\delta$ is a cut-off. In dimension $d=4$, we expect it to diverge as $\delta^p$ with $p$ depending on the Yang--Mills energy of the connection. It may be possible to define a renormalised version of this product in dimension $d=4$ and $d=5$.
\medskip

This formula is reminiscent of Symanzik's famous polymer representation \cite{Symanzik}, and loop-type representations have been used for partition functions in discrete settings early on, see \cite[Lemma 1.2]{BrydgesRW} but also the more recent \cite[Section 4]{Chandra}. Let us quickly mention some papers that we believe are the closest in spirit to our considerations. The body of work of Le Jan involves in the discrete setting analogues of many of the objects we consider and a wealth of results; for instance one may found a discrete analogue of \eqref{eq:detintro} in Section 6, Proposition 23 in \cite{lejan10} (see also Equation (11.1) in the upcoming book \cite{lejan22}). In \cite{KasselLevy}, from which the paper is partly inspired, the authors take this discrete formula further and prove Symanzik representations as we describe in the next paragraph. In \cite{LeJan}, results are obtained which are tantamount to ours for the case of a complex line bundle endowed with a flat connection. In the case of scalar Laplacians, the relation between their determinant and the Brownian loop soup is better understood, see e.g. \cite{Sheffield,Dubedat}. It is then the variation of the determinant with respect to the underlying metric $g$ which is considered, rather than the variation with respect to the connection $\nabla$. 

Our formula is innovative in the following ways. First, we consider a continuous space rather than a lattice or graph approximation. This turns a combinatorial problem into a analytical one; most of our work is precisely to deal with the issues inherent to the continuum. In the language of statistical mechanics, we \emph{start} in the ultraviolet limit.
Secondly, we consider all loops in the Brownian loop soup. We stress that the previous approaches often manage to discard loops of size less than $\delta$ in one way or another, for instance by imposing crossing events \cite{Dubedat} or non-contractibility \cite{LeJan}, or by working over a graph.
In cases where the small-$\delta$ limit is considered, we are only aware of results where a deterministic limit is taken \emph{after} averaging over the loop soup, for instance asymptotics for the renormalised total local time on average \cite{Sheffield}. On the other hand, although our product over the loop soup is also defined using a limiting procedure, it is ultimately a well-defined random variable which does not depend on any $\delta$. The method of construction and the analysis rely on regularisation, but the constructed object and results do not.
Most notably, the product over all loops gives a statistic-mechanical interpretation of these determinants.
Furthermore, our approach is quite general in terms of the geometric objects involved, in that we do not restrict to the Abelian setting (e.g. to line bundles), nor to topologically trivial of flat manifold, nor to trivial bundles.
\medskip

Let us break our proof into a few rough steps. We start by writing $\det(L)$ as an integral of a heat kernel. We then express this heat kernel in terms of holonomies along Brownian loops. Finally we recognise some space-time integrals involving a single Brownian loop as the expectation over the Brownian loop soup. This echoes ideas present in the literature from physics \cite{Simon}, index theory \cite{BGV,Gilkey,Rosenberg}, differential stochastic geometry \cite{Bismut,Norris,Hsu}, and conformal probability \cite{Werner,Kenyon}.

\medskip

Here are some properties of $\Det(L)$ that we can deduce, with little work, from the stochastic formula that we obtain. These are only a few examples, and we expect to exploit this formula further in future work.
\begin{itemize}
    \item If $E$ is a trivial bundle over $M$, a so-called \emph{diamagnetic inequality} states that the determinant of the Laplacian $\frac12\nabla^* \nabla+m$ is minimal when $\nabla$ is a trivial connection. From the $\zeta$-regularisation definition, this is a difficult result to obtain, and is for example the main result in \cite{Schrader} (also of \cite{Brydges} for the lattice case). Using Equation \eqref{eq:detintro}, this becomes trivial, noting that the traces must have modulus at most 1.
    \item We will also use this formula to show that $\Det(L_\nabla)$, as well as the heat kernel and Green kernel, depend continuously on $\nabla$ under the $\mathcal{C}^{2+\varepsilon}$ metric. We were not able to find this result in the geometric analysis literature, which makes us believe that a direct analytical approach would prove to be delicate. This is a cornerstone of our second result, namely a mathematically rigorous construction of a stochastic minimal coupling of a scalar field with respect to a random gauge field.
    \item As a last consequence, we prove in Section \ref{sec:inv} that ratios of determinants are invariant under conformal transformations in the case of surfaces.
    \item Although we do not work it out in this article, we believe that it is likely we can use this formula, together with some classical stochastic tools, to deduce properties of the asymptotic behaviour of the determinant of the magnetic Laplacian, associated with the connection $\d+i\alpha A$ on the trivial complex line bundle over $M$, in the limits $\alpha\to 0$ and $\alpha\to\infty$. Here the real-valued 1-form $A$ is a magnetic vector potential; see e.g. \cite{Raymond} for an introduction to magnetic Laplacians.
\end{itemize}

\paragraph*{Coupling scalar and gauge fields}
Motivated by physics, our second main goal is the following: given any probability measure $\mathbf{P}$ on smooth connections on $E$, define rigorously a measure on the couples $(\mathtt{section},\mathtt{ connection})$ over $M$, which informally is to be understood as
\begin{equation}
\label{eq:measureintro}
\frac{1}{Z} e^{-\frac14\|\nabla \phi\|^2-\frac12\langle \phi,m \phi\rangle}\mathcal{D}\phi \d \mathbf{P}(\nabla),
\end{equation}
where $Z$ is a normalising constant. As we will explain, the determinant of the Laplacian plays a crucial role in constructing such a measure. We will indeed define it rigorously, and give formula for some of its observables. With the tools developed at that point, it will be a relatively easy reap.

In physical terms, $\nabla$ corresponds to a gauge field (over the Riemannian space-time $M$), whilst $\phi$ corresponds to a scalar matter field.
The same expression with $\nabla$ replaced with the trivial connection $\d$ on a trivial bundle would lead to the interpretation of the marginal for $\phi$ as the measure of the Gaussian free field. The procedure of replacing $\d$ with $\nabla$ (and the trivial bundle with a generic one), i.e. ``replacing derivative with covariant derivatives,'' is what physicists loosely refer to as a \emph{minimal coupling}.

This measure is both a toy model and a building block for more involved and physically relevant couplings of the matter field $\phi$ with the gauge field $\nabla$.
In many contexts, one would like $\mathbf P$ to be the Yang--Mills measure on connections, which has support only on very irregular connections and falls way outside the scope of our framework, although it still allows one to consider for $\mathbf{P}$ mollifications (or otherwise regularised versions) of the Yang--Mills measure.
One would also like to consider potentials more general than the `kinetic energy' term ${\langle\phi, m\phi\rangle=\int m |\phi|^2}$; the Higgs field, for instance, involves a quartic term $\int_M|\phi|^4$, generally referred to as a \emph{self-interaction term.}

Notice that the equality \eqref{eq:detintro} has the advantage that it makes formal sense as a definition for the left-hand side, under the sole regularity assumption for $\nabla$ that holonomies along Brownian loops are well-defined objects. This is the case for some connections with a rather low regularity, for which it is very unclear that the $\zeta$-regularisation method is working. In fact, if the holonomies are also defined over Brownian bridges, we will show that the same is true for general polynomial observables of the measure \eqref{eq:measureintro}, in the sense of Theorem \ref{th:symanzik}. We do not construct the measure \eqref{eq:measureintro} under such low regularity conditions, but we do not exclude this should be possible.

We will show that the marginal $\mathbf{P}'$ of $\nabla$ under \eqref{eq:measureintro} is not the original probability measure $\mathbf{P}$, but $\mathbf{P}$ weighted by $\det(L)^{-\frac{1}{2}}$. In Liouville quantum gravity, it is known \cite{Sheffield} that weighting maps by powers of Laplacian determinants alters the central charge of the model, and thus produce a limit with a different regularity, which can also be done by changing the parameter $\gamma$ of the model, i.e. the exponent in the informal description $\colon e^{\gamma h}(\mathrm dx^2+\mathrm dy^2)\colon$.
Our stochastic representation suggests that the weight here tends to favour flatter connections, indeed suggesting that when $\mathbf{P}$ is a Yang--Mills field, the formally defined $\mathbf{P}'$ might correspond to a Yang--Mills field with a larger inverse temperature.
In \cite{Chandra}, it is shown, by controlling the partition functions in particular, that this weighted Yang--Mills measures, in discrete lattice, converge in the continuum limit at least along subsequences (and in the Abelian case).\footnote{The type of weight they consider is more general than a simple determinant of Laplacian, for they correspond to the more general $f(\phi)$ setting mentioned previously.} We believe it would be interesting to prove a similar result directly in the continuum, for example by also considering here mollifications of the Yang--Mills measure. For such an approach, and in particular to get estimates on the partition function, we think our formula, as well the bound given by Proposition \ref{prop:tech} below, would play a key role.

\paragraph*{Holonomy of small loops}
A key technical estimate for all our derivations is an asymptotic bound on the holonomy of a connection for Brownian loops in the $L^p$ norms. As $t$ goes to zero, the holonomy along a Brownian loop $\ell$ of duration $t$ should be close to the identity; our Proposition~\ref{prop:tech} asserts that the remainder $\mathcal Hol^\nabla(\ell)-\Id$ is of size $t$, and of size $t^2$ on average, in arbitrary dimension~$d$.
This estimation is of independent interest as it improves some heat kernel estimations. 

Some particular cases of this formula have a much simpler proof, for instance for a flat base manifold and abelian gauge group, or when the connection is flat. We present here the general case (for a principal bundle with compact fiber), which requires a non-trivial amount of work, and is not just about combining the different techniques to deal with these subcases. The control is obtained by probabilistic and mostly path-wise methods, rather than analytic ones.

\paragraph*{Outline}
Section \ref{sec:not} is devoted to the introduction of the objects we will consider throughout the paper, including some of their well-known properties. At the end of it, we formulate rigorously our first main result in the form of Theorem \ref{th:main1}, concerning the stochastic representation \eqref{eq:detintro} of the determinant.

Sections \ref{sec:rpz}, \ref{sec:mart} and \ref{sec:tech}  deal with the proof of this representation. The first two (Sections \ref{sec:rpz} and \ref{sec:mart}) are conditional to the estimate on the holonomies of small loops of Proposition \ref{prop:tech}, which we prove in Section \ref{sec:tech}.
In Section \ref{sec:rpz}, we prove that the difference of the $\zeta$-functions associated with two different connections over possibly different vector bundles $E$, the base manifold being fixed, can be written as a space-time integral involving holonomies of the connections over Brownian loops --- see Lemma \ref{le:zeta=E[chi]} for the precise result. In Section \ref{sec:mart}, we identify this integral as the expectation of a product over a Brownian loop soup under the form of Proposition \ref{prop:loops}. Together, they give Equation \eqref{eq:detintro} of the introduction, given in precise form in Theorem \ref{th:main1} of Section \ref{sec:not}.

We then turn to the construction of the measure described by \eqref{eq:measureintro}, which is the purpose of Section \ref{sec:Symanzik}. The measure is defined in Definition \ref{def:Symanzik}, making use of Proposition \ref{prop:measurability}. Theorem \ref{th:symanzik} gives a stochastic representation of polynomial functional of the associated scalar field, which is sometimes known as a Symanzik identity.

Our concluding Section \ref{sec:inv} describes the way the partition function changes under conformal transformations in dimension $2$, when the mass function is transformed appropriately:
\[ 
\frac{\det( L_{e^{2f}g,\nabla',e^{-2f}m})}{\det( L_{e^{2f}g,\nabla,e^{-2f}m})}=\frac{\det( L_{g,\nabla',m})}{\det( L_{g,\nabla,m}  )   }.
\]

\section{General framework}
\label{sec:not}

\subsection{Laplacians and their heat kernels, \texorpdfstring{$\zeta$}{zeta}-functions, and Green kernels}

In the following $M$ is a connected compact smooth manifold of dimension $d$, possibly with a smooth boundary, endowed with a Riemannian metric $g$ and the corresponding Levi-Civita connection $\nabla^{TM}$ on the tangent bundle $TM$. A \emph{mass} function $m\in \mathcal{C}^\infty(M,\mathbb{R}_+)$ is fixed. We say that $m$ vanishes if it is identically zero.
For integrals over $M$ we always write $\d x$ for the Riemannian volume element $\d \vol_g(x)$. 
We consider $E$ a real or complex vector bundle over $M$, with finite rank $n$ and endowed with a metric, either Euclidean or Hermitian.

A metric connection $\nabla$ on $E$ gives, for every smooth curve $\gamma:[0,T]\to M$, 
a notion of parallel transport along it, that we denote by $\mathcal Hol^{\nabla}(\gamma)$. It is an isomorphism from $E_{\gamma(0)}$ to $E_{\gamma(T)}$, and in fact an isometry since $\nabla$ is metric. For $0\leq s\leq t\leq T$, we also write $\mathcal Hol^{\nabla}_{s,t} (\gamma)$ for $\mathcal Hol^{\nabla} (\gamma_{|[s,t]})$, and $\mathcal Hol^{\nabla}_{t,s} (\gamma)$ for $\mathcal Hol^{\nabla}_{s,t} (\gamma)^{-1}$, which is also equal to $\mathcal Hol^{\nabla}_{s,t} (\gamma)^*$. In particular, for any $\gamma$ the curve $t \mapsto \mathcal Hol^{\nabla}_{0,t} (\gamma)$ is the $\nabla$-horizontal lift of $\gamma$ to $E$, that is $\mathcal Hol^{\nabla}_{0,t} (\gamma)$ is the isomorphism from $E_{\gamma(0)}$ to $E_{\gamma(t)}$ such that $\mathcal Hol^{\nabla}_{0,t} (\gamma)(v)$ is the parallel transport of $v$ along $\gamma$.
Generalities regarding principal bundles and parallel transport may be found in \cite[Chapter 8]{Spivak2}. For their stochastic counterparts when the driving curve is a diffusion process, see \cite[Section 2.3]{Hsu}.

For $x\in M$ and $u\in \End(E_x)$, we define $\tr(u)$ as the normalised trace of $u$: for any orthonormal basis $(e_i)$ of $E_x$,  $\tr(u)=\frac{1}{n}\sum_{i=1}^n\langle u(e_i),e_i\rangle $.

If $M$ does have a boundary, we think of it as a domain with smooth boundary in some closed manifold $\widehat M$, for instance by considering the double $M\times\{0,1\}/\{(x,0)\sim(x,1),x\in\partial M\}$ of $M$ together with its smooth structure, see e.g. \cite[Example 9.32]{Lee}. In that case, we can assume without loss of generality that all the objects defined above are actually restrictions of objects defined on $\widehat M$. We say that a smooth section of some vector bundle over $M$ vanishes at the boundary if it is the restriction of a smooth section over $\widehat M$ that vanishes on $\partial M$ (up to an isomorphism identifying the bundle on $M$ with the restriction of a bundle on~$\widehat M$). The set of such sections is written $\Gamma_0(M,E)$.

All the metrics and connections that we consider through the paper are smooth, and we always assume compatibility, i.e. every connection is metric. Here, as everywhere above and below, smooth means infinitely differentiable.
\bigskip

Associated with $\nabla$ (and $m$)
is a Dirichlet Laplace-type operator $L_{\nabla}$ acting on smooth sections of $E$ vanishing on the boundary, defined by
\begin{equation}
\label{eq:Ldef}
L_{\nabla}s\coloneqq -\frac{1}{2}\Tr(\nabla^{E\otimes T^*M}\nabla)(s)+ ms,
\end{equation}
where the connection $\nabla^{E\otimes T^*M}$ on the vector bundle $E\otimes_M T^*M$ over $M$ is defined as 
$\nabla^{E\otimes T^*M}\coloneqq \nabla\otimes\operatorname{Id}_{T^*M} +\operatorname{Id}_E\otimes\nabla^{T^*M}$, where $\nabla^{T^*M}$ is the connection on $T^*M$ induced from $\nabla^{TM}$ by duality.
For $s$ compactly supported inside $\operatorname{int}M$, $Ls$ is also given by $Ls=\frac{1}{2}\nabla^*\nabla s +ms$, where $\nabla^*$ is the formal adjoint of $\nabla$: for any two smooth sections $s,s'$ compactly supported in the interior of $M$,
\[ \langle L s, s'\rangle_{L^2(M,E)}
 = \frac{1}{2} \langle \nabla s , \nabla s' \rangle_{L^2(M,T^*M\otimes E)}+ \langle m s, s'\rangle_{L^2(M,E)}. \]
For the expression of $L_\nabla$ in coordinates, see e.g. \cite[Formula (4.1.15)]{Gilkey}.

Sometimes we compare the operator associated with different connections (or even different bundles), in which case it is convenient to keep this index $\nabla$, but when we consider a single connection $\nabla$, we simply drop this index (or superscript) $\nabla$, both for $L_\nabla$ and for the other notations we will define.

For the spectral theoretic aspects of our analysis, we consider the Dirichlet extension of $L$. Namely, starting from the operator acting on sections vanishing on the boundary, we also write $L$ (or $L_\nabla$) for its unique self-adjoint extension, whose domain is $H^2(M,E)\cap H^1_0(M,E)$. It has a discrete spectral resolution $(\phi_i,\lambda_i)_{i\in \mathbb{N}}$ with smooth, $L^2$--orthogonal eigensections $\phi_i$, and non-negative eigenvalues $\lambda_i$ whose only accumulation point is at $+\infty$ (see e.g. \cite[Lemma 1.6.3]{Gilkey} or \cite[Section 8.2]{TaylorII}).

We will often need the semigroup to be exponentially contracting, i.e. we want all the eigenvalues of $L$ to be larger than a positive constant. Because of the considerations above, it is equivalent to say that the kernel is trivial. In Remarks \ref{rk:kernel} and \ref{rk:kernel2}, we discuss geometric constraints enforcing this condition; for now, let us mention that the mass being positive anywhere, or the boundary being non-empty, are both sufficient conditions. These are the most important ones, because \emph{all} semigroups becomes exponentially contracting, most notably the one corresponding to the Laplace--Beltrami operator acting on functions, which means that almost surely, the massive Brownian loop soup has a finite number of large loops (see below).

There is a heat equation associated to $L$, and we denote by $(e^{-tL})_{t\geq0}$ its semigroup, seen for instance as a collection of compact operators over $L^2(E)$. For all $t>0$, the operator $e^{-tL}$ admits (see e.g. \cite{Greiner}, or \cite[Lemma 1.6.5]{Gilkey} for a simpler proof in the case without boundary) a smooth kernel
\[K^\nabla_t\in \Gamma(M\times M, E\boxtimes E^*)\] 
defined by the equality
\begin{equation}
\label{eq:Kdef}
K_t(x,y)= \sum_i e^{-t \lambda_i} \phi_i(x)\otimes \phi_i^*(y), 
\end{equation}
where the sum is absolutely convergent in all $\mathcal C^k([\varepsilon,\infty)\times M\times M,E\boxtimes E^*)$.

Identifying $E_x\otimes E_y^*$ with the space of linear operators $E_y\to E_x$ via $(e\otimes \phi)(v)=\phi(v)e$, we can integrate any regular section $s$ of $E$ against the kernel $K$, and we get the following expression for the semigroup generated by $-L$:
\[ (e^{-t L } s)(x)=  \int_M K_t(x,y)(s(y)) \d x.
\]

When $M$ has no boundary, as $t\to 0$, the heat kernel admits the asymptotic expansion 
\begin{equation}\label{eq:Kexpansion}
K_t(x,x)= \sum_{i=0}^{k} u_i(x,x) t^{i-d/2} +O( t^{k+1-d/2}), 
\end{equation}
where $k$ is an arbitrary non-negative integer, the sections $u_i$ are smooth, $u_0$ is nowhere vanishing, and the remainder is uniform in $x\in M$, see e.g. \cite[Section 3.3]{Rosenberg}. When $M$ has a non-empty boundary, this expansion fails from being uniform near the boundary.

The trace $\Tr(e^{-tL})=\sum_i e^{-t\lambda_i}$ is finite and given by 
\[ \Tr(e^{-tL})= n \int_M \tr( K_t(x,x))\d x.\]
Indeed, 
\begin{align*}
n \int_M \tr( K_t(x,x))\d x
&= \int_M \sum_i e^{-t \lambda_i}  n \tr(  \phi_i(x)\otimes \phi_i^*(x)  )\d x\\
&= \sum_i e^{-t \lambda_i} \int_M  \|\phi_i(x)\|^2\d x=\sum_i e^{-t \lambda_i}.
\end{align*}
The asymptotic expansion \eqref{eq:Kexpansion} can be integrated over $M$, and in fact even in the case when $M$ has a non-empty boundary, there exists constants $a_i$ such that for all $k\geq0$, we have
\begin{equation}
\label{eq:KexpansionIntegrated}
\Tr(e^{-tL})= \sum_{i=0}^k a_i t^{(i-d)/2}+O(t^{(k+1-d)/2}).
\end{equation}
If $M$ has no boundary, $a_{2i}=0$ for all $i$. See e.g. \cite[Theorem 2.6.1]{Greiner}.
 
We write $\pi_L$ for the orthonormal projection from $\Gamma(M,E)$ to $\Ker(L)$. The Weyl asymptotics for the eigenvalues of $L$ (which can be deduced from the main term of \eqref{eq:KexpansionIntegrated}, following e.g. \cite[Section 15]{Shubin} or \cite[Lemma 2.43]{BGV}) ensures that the $k^{\text{th}}$ eigenvalue grows as $k^{2/d}$. Thus, for
$\mathfrak{R}(s)> \frac{d}{2}$, the sum
\[ 
\zeta_\nabla(s)\coloneqq \sum_{i: \lambda_i\neq 0} \lambda_i^{-s},
\]
is absolutely convergent. Using the identity
\[ \lambda^{-s}= \frac{1}{\Gamma(s)}\int_0^\infty  e^{-t\lambda} t^{-1+s}\d t,  \]
 valid for any $\lambda$ and $s$ with $\lambda>0$ and $\mathfrak{R}(s)>0$, it rewrites as
\[
\zeta_\nabla(s)=\frac{1}{\Gamma(s)}\int_0^\infty  \Tr(e^{-tL}-\pi_L) t^{-1+s}\d t
\]
in the same range $\mathfrak{R}(s)> \frac{d}{2}$. The integral and sum can be exchanged since the real case $s>d/2$ gives a summable upper bound.

This function defined on the half-plane $\mathfrak{R}(s)> \frac{d}{2}$ admits a meromorphic extension to the whole plane, whose poles are contained in
$ \{\frac{d}{2}, \frac{d-1}{2}, \dots, \frac{1}{2} \}\cup \{-\frac{1}{2},-\frac{3}{2},\dots,\}$ (the proof of \cite[Theorem 5.2]{Rosenberg} in the case $\partial M=\emptyset$ translates readily to our situation using \eqref{eq:KexpansionIntegrated}%
\footnote{This is in fact a general statement: for any continuous function $f=f(t)$ defined on $(0,\infty)$, exponentially decreasing as $t\to + \infty$, and for which there exists an increasing sequence of real numbers $\beta_i$ which converges toward $+\infty$  and a sequence $a_i$ of non-zero complex numbers such that for all $k$, $f(t)=\sum_{i=1}^k a_i t^{\beta_i} +O(t^{\beta_{k+1}})$  as $t\to 0$, the Mellin transform  of $f$ is meromorphic in the plane with poles exactly at the $\beta_i$, which are all simple poles. This is proven by considering a half-plane $\Re(s)>-b$, and  splitting the integral which defines the Mellin transform at $t=1$. The exponential decay ensures that the large $t$ integral is holomorphic on the plane. For the small $t$ integral,
each term $a_it^{\beta_i}$ contributes a simple pole at $-\beta_i$, and when $b<\beta_{k+1}$ the remainder is holomorphic on the considered half-plane. Remark that for the $\zeta$ function, possible poles at negative integer values are cancelled by the pole of $\Gamma$.}\textsuperscript{,}%
\footnote{ We think there might be a mistake in \cite[Theorem 5.2]{Rosenberg}. We do not see why there would be no poles at negative half-integers. In any case, we do believe these poles exist in the case with boundary. Remark in this paper we only need to know about the behavior of $\zeta$ in the half plane $\Re(z)>-\varepsilon$ anyway.}).

Provided the first eigenvalue is positive, the Green kernel
\[ G_\nabla\in \Gamma(M\times M,E \boxtimes E^*) \]
is defined outside the diagonal $\{(x,x), x\in M\}$ by \[G_\nabla(x,y)\coloneqq\int_0^\infty K_t(x,y) \d t.\]
The operator of convolution with $G_\nabla$ is an inverse of $L_\nabla$, in the sense that for any smooth section $s\in \Gamma(M,E)$, it holds that
\[
s= L_\nabla( G_\nabla s ),\qquad (G_\nabla s) (x)\coloneqq \int_{M} G_\nabla(x,y)s(y)\d y.
\]

We call $p_t$ the heat kernel for half of the massless Laplace--Beltrami operator $\Delta$ on $\widehat M$, which acts on real-valued functions over $M$. It corresponds to the case $E=\widehat M\times \mathbb{R}$, $m\equiv0$ and $\nabla^E=\d$, so it satisfies all the previous estimates (with rank $n=1$ for the bundle). Furthermore the kernel of the Laplace--Beltrami operator is generated by a single constant function, so $p_t(x,x)$ converges exponentially fast towards $\frac{1}{\Vol_g(M)}$ as $t\to + \infty$, uniformly in $x\in M$. As for the small times, from the expansion \eqref{eq:Kexpansion} we get the useful control
\begin{equation}\label{eq:unifpbound}
\sup_{t\leq1}\sup_{x\in M}t^{-d/2}p_t(x,x)<\infty.
\end{equation}

\subsection{Brownian bridge and loop soup}
In this paragraph, we recall the (slightly amended) definition of the loop soup introduced in \cite{Werner}. We use the notation $\mathbb{E}_{t,x,y}$ for the expectation with respect to $\mathbb{P}_{t,x,y}$, under which $W$ is a Brownian bridge in $\widehat M$ from $x$ to $y$ with duration $t$, and $\mathbb{E}_x$ (resp. $\mathbb E_{t,x}$) for the expectation with respect to $\mathbb{P}_{x}$ (resp. $\mathbb P_{t,x}$), under which $W$ is a Brownian motion in $M$ from $x$ (resp. with duration $t$). One can think of $W$ under $\mathbb{P}_{t,x,y}$ as the solution to some stochastic differential equation driven by a Brownian motion, see Appendix \ref{app:Econt}. For a function $m\in \mathcal{C}^\infty(M,\mathbb{R}_+)$, we define $\mathbb{P}^m_{t,x,y}$ the measure (with total mass not necessarily equal to $1$) given by
\[
\mathbb{P}^m_{t,x,y}(A) \coloneqq \int_A \mathbf{1}_{\text{$W$ stays in $M$}} \cdot e^{-\int_0^t m(W_s) \d s }\d \mathbb{P}_{t,x,y}(W).
\]
It can be thought of as the same expectation with the indicator removed, extending $m$ by infinity outside of $M$. We write $\mathbb{E}^m_{t,x,y}$ for the corresponding integral, which we treat as an expectation:
\[ 
\mathbb{E}^m_{t,x,y}[f(W)]\coloneqq \int_{\mathcal{C}([0,t],M)}   f(W) \d \mathbb{P}^m_{t,x,y}(W)=\int \mathbf{1}_{\text{$W$ stays in $M$}} \cdot e^{-\int_0^t m(W_s) \d s } f(W) \d \mathbb{P}_{t,x,y}(W).
\]

A loop on $M$ is an element $\ell$ of the set
\[
\mathscr{L}\coloneqq \bigcup_{T>0} \{ \ell\in \mathcal{C}([0,T],M ) : \ell(0)=\ell(T) \}.
\]
For a loop $\ell$, we write $T(\ell)$ for the unique $T$ such that $\ell \in\mathcal{C}([0,T],M)$ and we call it the duration of $\ell$.

There are many natural topologies on $\mathscr L$, for instance the naive topology coming from the bijection
\[ (T,\gamma)\in(0,\infty)\times\mathcal C(\mathbb S^1,M)\mapsto \big(T,\gamma(\cdot/T)\big),  \]
or more subtle topologies that forget about the basepoint (but not the orientation!) and time parameter. We will need very little from the topology (for instance we only care about the $\sigma$-algebra), so we leave it to the readers to choose their own. We ask however that the above map be measurable with respect to the product Borel algebra on the left, and that integrals of smooth functions and forms be measurable in the following sense. For every (smooth) function $f:M\to\mathbb R$, and every connection $\nabla$ on every $E$, we ask that some measurable maps
\begin{align*} \int_\bullet f:\mathscr L\to\mathbb R, && \mathcal Hol^\nabla:\mathscr L&\to\mathrm{End}_M(E) \end{align*}
exist, that coincide almost surely under (the pushforward to $\mathscr L$ of) all $\mathbb P^m_{t,x,x}$ with, respectively, the integral of $f$ along the trajectory, and its Stratonovich holonomy. Note that the duration of a curve is simply the integral of the constant function $\mathbf1$, hence measurable in our sense. By the usual theory of stochastic differential equations on manifolds, the naive topology is such an option.

For a mass function $m\in \mathcal{C}^\infty(M,\mathbb{R}_+)$, the massive Brownian loop soup measure $\Lambda$ with mass $m$ (and with intensity $1$) is the measure on $\mathscr{L}$ given by 
 \[
 \Lambda(A)= \int_0^\infty  \int_M \frac{p_t(x,x) }{t} \mathbb{P}^m_{t,x,x}(A) \d x\d t.
 \]
The Brownian loop soup $\mathcal{L}$ on $M$ with mass $m$ and intensity $\alpha$ is by definition a Poisson point process on $\mathscr{L}$ with intensity $\alpha\Lambda$. Technically speaking, the version of \cite{Werner} is concerned with the case $m=0$; the mass term is natural and has been considered e.g. in \cite{Camia} with the convention $m_\mathrm{PS}=m_\mathrm{Camia}^2$. For notational simplicity, we will assume $\alpha=1$ in most of the paper, since there is no added difficulty at all to deal with the general case.  We fix a probability space endowed with a Brownian loop soup, and write $\mathbb{P}$ and $\mathbb{E}$ for the corresponding measure and expectation, or $\mathbb{P}^{\mathscr L}_\alpha$, $\mathbb{E}^{\mathscr L}_\alpha$ in the rare cases when we want to explicitly write the dependency in $\alpha$.

For $0<\delta< R< \infty$, we define $\mathscr{L}_\delta$ (resp. $\mathscr{L}^R$, resp. $\mathscr{L}_\delta^R$ as the subset of $\mathcal{L}$ of loops with duration at least $\delta$ (resp. less than $R$, resp. in between $\delta$ and $R$). 
Similarly we define the measures $\Lambda_\delta$, $\Lambda^R$, and $\Lambda_\delta^R$ as the restrictions of $\Lambda$ to the corresponding subspaces of $\mathscr L$, and the point processes $\mathcal{L}_\delta$, $\mathcal{L}^R$, and $\mathcal{L}_\delta^R$ as the intersections of $\mathcal{L}$ with the corresponding subsets of $\mathscr L$. The latter point processes are Poisson processes with intensities given by the former measures.

Notice the measure $\Lambda$ has infinite mass whilst $\Lambda_\delta^R$ has finite mass for all $0<\delta<R<+\infty$. Thus, the Poisson process $\mathcal{L}$ almost surely contains infinitely many loops, but $\mathcal{L}_\delta^R$ is almost surely finite for all $\delta>0, R<\infty$. In the case when $m$ is not identically vanishing or $M$ has a boundary, $\mathcal{L}_{\delta}$ is also finite.

\subsection{Stochastic representations of determinants}

Our representation result for determinants of Laplacians starts with comparing two zeta functions. For this purpose, let us consider two metric connections $\nabla_0$ and $\nabla_1$ on two metric bundles $(E_0,h_{E_0})$ and $(E_1,h_{E_1})$. These bundles possibly have different ranks $n_0$ and $n_1$. To these objects correspond two Laplace operators $L_0$ and $L_1$ defined as above. For a Brownian loop $W$ with duration $t$, set 
\begin{equation}\label{eq:defXi}
\chi(W)=\chi_{\nabla_0,\nabla_1}(W)\coloneqq  \tr (\mathcal Hol^{\nabla_1}(W) )-\tr (\mathcal Hol^{\nabla_0}(W) ).
\end{equation}

\begin{theorem}
\label{th:main1}
Suppose that $M$ has dimension $d\in\{2,3\}$.
Assume that
\begin{equation}
\tag{$\ast$} \label{eq:kernels}
\frac{1}{n_0}\dim( \Ker L_0)=\frac{1}{n_1}\dim( \Ker L_1). \end{equation}
Then
\begin{equation}
\label{eq:zetaprimeisadoubleintegral}
\frac{1}{n_1}\zeta'_1(0)-\frac{1}{n_0}\zeta'_0(0)= \int_0^\infty \int_M \mathbb{E}^m_{t,x,x}[\chi(W)]  \frac{p_t(x,x)}{t}\d t \d x.
\end{equation}

Make the stronger assumption that $d\in\{2,3\}$ and that either $m$ is not constant equal to zero or $\partial M$ is non-empty.
Then $\mathcal{L}_{\delta}$ is almost surely finite for all $\delta>0$, and the limit
\[ \prod_{\ell\in \mathcal{L}} \big(1+  \chi(\ell)\big) \coloneqq \lim_{\delta\to0}\prod_{\ell\in \mathcal{L}_{\delta}} \big(1+  \chi(\ell)\big)  \]
of finite products exists in the almost sure sense and in the $L^p(\mathbb P^\mathscr{L}_\alpha)$ sense for all $p\geq1$ and $\alpha>0$.
For every intensity $\alpha>0$, it satisfies
\begin{equation}
\label{eq:pre-ratio}
\exp\Big(\frac{\alpha}{n_1}\zeta'_1(0)-\frac\alpha{n_0}\zeta'_0(0) \Big)
    =\exp\Big(\alpha\int_0^\infty \hspace{-0.2cm} \int_M \mathbb{E}^m_{t,x,x}[ \chi(W)   ]  \frac{p_t(x,x)}{t}\d x \d t \Big)
    =\mathbb{E}^\mathscr{L}_\alpha\bigg[\prod_{\ell\in \mathcal{L}}\big(1 + \chi(\ell)\big) \bigg].
\end{equation}
\end{theorem}

\begin{remark}
It is part of the result that the right-hand side in \eqref{eq:zetaprimeisadoubleintegral} is well-defined. In fact, $\mathbb{E}^m_{t,x,x}[| \chi(W)|]<\infty$ for all $t>0$ and $x\in M$, and 
$\int_0^\infty \int_M |\mathbb{E}^m_{t,x,x}[ \chi(W)]| \frac{p_{t}(x,x)}{t}\d t \d x<\infty$.
The integral features the intensity of the loop soup, and \emph{formally} \eqref{eq:zetaprimeisadoubleintegral}${}=\alpha\int_{\mathscr{L}} \chi(\ell) \Lambda(\mathrm{d} \ell)$, but this last expression is in fact ill-defined as $\int_{\mathscr{L}} |\chi(\ell)| \Lambda(\mathrm{d} \ell)=+\infty$. However, we do have
\begin{align*}
   \int_0^\infty \int_M \mathbb{E}^m_{t,x,x}[\chi(W)]  \frac{p_t(x,x)}{t}\d t \d x 
&= \lim_{\delta\to0}\int_\delta^\infty \int_M \mathbb{E}^m_{t,x,x}[\chi(W)]  \frac{p_t(x,x)}{t}\d t \d x\\
&= \lim_{\delta\to 0} \int_{\mathscr{L}_\delta} \chi(\ell) \Lambda(\mathrm{d}\ell).
\end{align*}
\end{remark}

When $E_0$ is the trivial line bundle over $M$, the theorem implies the following.
\begin{corollary}
\label{cor:main}
    Assume that either $m\neq 0$ or $\partial M\neq \emptyset$. Let $E$ be a vector bundle over $M$ with rank $n$, endowed with a bundle metric, and $\nabla$ be a metric connection over $E$. Let $\zeta$ be the $\zeta$-function associated with $L=\frac{1}{2}\nabla^*\nabla+m$ in the sense of equation \eqref{eq:Ldef}. Then the limit
    \[ \mathcal{Z} = \prod_{\ell\in \mathcal{L}} \tr(\mathcal Hol^\nabla(\ell))
    \coloneqq \lim_{\delta\to0}\prod_{\ell\in \mathcal{L}_\delta} \tr(\mathcal Hol^\nabla(\ell))\]
    of finite products converges in the almost sure sense and in the $L^p(\mathbb P^\mathscr{L}_\alpha)$ sense for all $p\geq 1$ and $\alpha>0$. Furthermore,
    \[ \exp\Big( \frac{\alpha}{n} \zeta'(0) \Big)= (C_{g,m})^\alpha\,\mathbb E^\mathscr{L}_\alpha\big[\mathcal{Z} \big],   \]
    where $C_{g,m}=\exp( \zeta'_0(0))$ for $\zeta_0$ the $\zeta$-function associated with $\Delta/2+m$, where $\Delta$ is the Dirichlet Laplace--Beltrami operator on $M$. 
\end{corollary}

This is the expected formula \eqref{eq:detintro} of the introduction.

\begin{remark}
It is fairly common (see e.g. \cite{Bismut,Hsu,BGV}) to consider Laplace-type operators for which the mass term is not a scalar function but a section of $\End(E)$, usually involving a curvature (and described by a Weitzenböck formula, see e.g. \cite[Theorem 2.38]{Rosenberg}). This is the case for example when one considers the Hodge Laplacian, a conformal Laplacian (acting on an hermitian bundle), or a spin Laplacian.
In fact, any symmetric differential operator with smooth coefficients and principal symbol $\frac12g^{-1}\otimes\operatorname{id}_E$ is equal to some $L_{\nabla,m}$, where $m$ must be allowed to be endomorphism-valued, see \cite[Lemma 4.1.1]{Gilkey}. 

In that case, it is not possible to keep the probabilistic interpretation of the mass as a killing rate weighting down the intensity of the loop soup, and we must consider the full loop soup. Instead, the mass term has to be taken into account in the holonomy terms, as an extra drift term in the SDE defining the holonomy; the corresponding quantity in the discrete setting is referred to as the \emph{twisted holonomy} by \cite{KasselLevy}. 
Our analysis should extend to this case, under the assumption that the eigenvalues of the endomorphism are pointwise non-negative. We expect Theorem \ref{th:main1} to still hold with this modification of $\mathbb E^m$, $\mathcal L$ and $\mathcal Hol^\nabla$, but only for $E_0=E_1$.
\end{remark}

\section{Heat kernel representation of the determinant}
\label{sec:rpz} 

In this section, we prove the first formula in Theorem \ref{th:main1}, conditionally on the following asymptotic estimation, which will be used in this section and the next, but will only be proved in section \ref{sec:tech}. All of these inequalities are actually easy consequences of the first.
\begin{prop}
\label{prop:tech}
There exists $C>0$ which depends on $E$ and $\nabla$ such that for all $t>0$ and $p\geq 1$, 
\[ \sup_{x\in M}  \mathbb{E}_{t,x,x}[|\operatorname{Id}_{E_x}-\mathcal{H}ol^{\nabla}_{0,t}(W)|^p ] \leq (Cpt)^p.
\]
Furthermore, for all $t>0$, 
\[ \sup_{x\in M}  \big|\mathbb{E}_{t,x,x}[\operatorname{Id}_{E_x}-\mathcal{H}ol^{\nabla}_{0,t}(W) ]\big| \leq Ct^2.
\]
Identical bounds hold when $\mathbb{E}_{t,x,x}$ is replaced with $\mathbb{E}_{t,x,x}^m$, and
similar bounds hold for $\chi_{\nabla_0,\nabla_1}$: there exists $C$ which depends on $E_1,E_2,\nabla_1,\nabla_2$ such that 
\begin{align*}
\sup_{x\in M}\mathbb{E}^m_{t,x,x}\big[  \big|\chi_{\nabla_0,\nabla_1}(W)\big|^p\big]&\leq (Cpt) ^p
&& \text{and} &
\sup_{x\in M}  \big|\mathbb{E}^m_{t,x,x}\big[\chi_{\nabla_0,\nabla_1}(W)\big]\big|&\leq C t^2.
\end{align*}
The constant now depends on both connection, but not on $m$.
\end{prop}
    
\begin{remark}
For comparison, let us remark that in local coordinates such that $\nabla=\d+A$, and with $W$ a Brownian motion with duration $t$, the quantity $|\mathcal{H}ol^{\nabla}_{0,t}(W)-\Id_{E_x}|$ is typically of order $\sqrt{t}$ rather than $t$, even in the case when $M$ is flat and $E$ is the trivial complex line bundle. It is only because we are taking the holonomy along a Brownian \emph{loop}, rather than a Brownian \emph{path} with different endpoints, that we can save an extra factor $\sqrt{t}$. 
    
Let us consider the trivial case $E=\mathbb{R}^2\times \mathbb{C}$, to understand why we can indeed save an extra factor. In this case, the endomorphism $A$ can naturally be identified with a real-valued $1$-form, and $\mathcal{H}ol^{\nabla}_{0,t}(W)$ is given by $\exp(i\int_0^t A_{W_s} \d W_s)$, where the integral is to be understood in the sense of Stratonovich (the factor $i$ comes from the identification of $\mathfrak{u}(1)$ with $i\mathbb{R})$. When $W$ is a loop, one can replace $A$ with $A'=A+\grad f$, for an arbitrary smooth function $f$, without changing the value of the integral. In particular, one can choose $f$ such that, at the point $W_0$, $\grad f=-A$. Since $A$ (and thus also $A'$) is Lipschitz continuous and $|W_s-W_0|$ is at the most of order $\sqrt{t}$, we deduce that $A'_{W_s}$ will be typically of order $\sqrt{t}$, whilst $A_{W_s}$ is typically of order 1. This replacement of $A$ with $A'$ is of course not possible when we consider a Brownian path that is not a loop. 

This simple idea is in fact the starting point for the proof of Proposition \ref{prop:tech}, not only for the trivial case but for the general one.  
\end{remark}

As an entry point, we will use the following Feynman-Kac type formula:
\begin{theorem}
    \label{th:K=pE}
    For all $x,y\in\operatorname{int}M$, $t>0$, the kernel $K_t$ defined by \eqref{eq:Kdef} decomposes as
    \[K_t(x,y)=p_t(x,y)\mathbb E^m_{t,x,y}\big[ \mathcal Hol^\nabla_{t,0}(W)\big].
    \]
\end{theorem}

The case where $M$ has no boundary and the equality holds only for almost every $y$ can be found in the work of Norris \cite[equation (34)]{Norris}. The case with boundary is essentially the same, up to a stopping time argument that we present here for convenience. 

\begin{proof}
The expression \eqref{eq:Kdef} is uniformly convergent over positive times, so we know that $K$ and $p$ are continuous (in fact they are smooth). We prove in Lemma \ref{le:Econt} of Appendix \ref{app:Econt} that the expectation on the right hand side is continuous with respect to the triple $(t,x,y)$, hence we only have to prove the equality almost everywhere.

We follow the outline of proof of \cite[Equation (34)]{Norris}. As discussed in Section \ref{sec:not}, we can assume that all our objects (e.g. $E$ and $\nabla$) are restrictions to $M$ of smooth objects defined over $\widehat M$ (e.g. $\widehat E$ and $\widehat\nabla$). Corresponding to all these parameters is a Laplace operator $\widehat L=\frac12\widehat\nabla^*\widehat\nabla+\widehat m$ defined on smooth sections of $\widehat E$ (which may not be non-negative since we may not be able to ensure non-negativity of $\widehat m$).

Let $\hat h:(t,x)\mapsto\hat h(t,x)$ be a time-dependent section of $\widehat E$. Define the process
\[Z^{\hat h}:t \mapsto \exp\Big(-\int_0^t m(X_s)\d s\Big) \mathcal Hol^{\widehat\nabla}_{t,0}(X)\hat h(t,X_t)\]
with values in $\widehat E_x$, where $X$ is a Brownian motion in $\widehat{M}$ starting from $x\in M$. If $(t,x)\mapsto\hat h(t,x)$ is smooth, then by It\^o's formula (see the details in \cite{Norris}), $Z^{\hat h}$ is a semimartingale satisfying
\[ \d Z^{\hat h}_t = \exp\bigg(-\int_0^t m(X_s)\d s\bigg) \mathcal Hol^{\widehat\nabla}_{t,0}(X) ( \widehat L\hat h-\partial_t\hat h)(t,X_{t})\d t +\d Y^{\hat h}_t,  \]
where $Y^{\hat h}$ is a local martingale. For a fixed time horizon $T$, we see that $(Y^{\hat h}_t)_{t\leq T}$ is bounded, so it is actually a martingale.

Let $s\in\mathcal{C}^\infty_0(M,E)$. Let $\tau$ be the first time when $X$ hits $\partial D$, and set
\[ h:(t,x)\mapsto \big(e^{-(T-t)L}s\big)(x) = \int_MK_{T-t}(x,y)s(y)\mathrm dy \]
for $x\in M$; in particular $h(t,x)=0$ for $x\in\partial M$. Let $(C_n)_{n\geq0}$ be a compact exhaustion of $\operatorname{int}M$ with $x\in\operatorname{int}C_0$, $\tau_n$ the exit time of $C_n$ and $h_n$ a smooth time-dependent section of $\widehat E$ that coincides with $h$ over $(-\infty,T-1/n)\times C_{n+1}$. Applying the above, $t\mapsto Z^{h_n}_{t\wedge\tau_n}$ is a martingale, and we get
\[
  \mathbb{E}[Z^h_{T\wedge\tau}]
= \lim_{n\to\infty}\mathbb{E}[Z^{h_n}_{T\wedge\tau_n}]
= \lim_{n\to\infty}\mathbb{E}[Z^{h_n}_0]
= h(0,x)
= \big(e^{-TL}s\big)(x).
\]
On the other hand,
\begin{align*}
\mathbb{E}[Z^h_{\tau\wedge T}]
&= \mathbb{E}[\mathbf{1}_{T<\tau} Z^h_{T}]+0\\
&= \int_{M} p_T(x,y) \mathbb{E}_{T,x,y}[\mathbf{1}_{T<\tau} e^{-\int_0^T m(X_t) \d t }\mathcal Hol^\nabla_{T,0}(X)( s(y) )  ] \d y \\
&=    \int_{M} p_{T}(x,y) \mathbb{E}^m_{T,x,y}[\mathcal Hol^\nabla_{T,0}(X)  ]s(y) \d y .
\end{align*}
Thus, for all $x\in M$ and $T>0$, for all $s\in \mathcal{C}^\infty_0(M)$,
\[ \big(e^{-TL}s\big)(x)=\int_{M} p_{T}(x,y) \mathbb{E}^m_{T,x,y}[\mathcal Hol^\nabla_{T,0}(X) ] s(y)\d y.  \]
Since $K_T$ is the kernel of $e^{-TL}$, we arrive to the equality
\[ K_T(x,y) = p_T(x,y)\mathbb E^m_{T,x,y}[\mathcal Hol^\nabla_{T,0}(X)] \]
for all $T>0$ and $x\in\operatorname{int}M$, and almost all $y\in\operatorname{int}M$.
\end{proof}

\begin{corollary}
\label{coro:integrability}
For all $s>d/2-2$, the integral
\[ \int_0^1\int_M\Big|\operatorname{tr}(K_t(x,x))-p_t(x,x)\mathbb E^m_{t,x,x}[\mathbf1]\Big|\mathrm dx \frac{\mathrm dt}{t^{1-s}} \]
is finite --- note that $\mathbb E^m_{t,x,x}[\mathbf1]$ is simply the total mass of $\mathbb P^m_{t,x,x}$.\medskip

For all $s\in\mathbb R$, the integral
\[ \int_1^\infty\int_M\bigg|\operatorname{tr}K_t(x,x)-\frac{\dim\big(\ker L\big)}{n\cdot\operatorname{vol}(M)}\bigg|\mathrm dx\frac{\mathrm dt}{t^{1-s}}\]
is finite.
\end{corollary}

\begin{proof}
For $t$ small, according to Theorem \ref{th:K=pE},
the integral is
\[ \int_0^1\frac{\mathrm dt}{t^{1-s}}\int_M\big|\operatorname{tr}\mathbb E^m_{t,x,x}\big[\mathcal Hol^\nabla_{0,t}(W)-\mathrm{Id}_{E_x}\big]\big|\cdot p_t(x,x)\mathrm dx. \]
The trace is of order $t^2$ by Proposition \ref{prop:tech}. As for the heat kernel, we know from equation \eqref{eq:unifpbound} that it must be of order $t^{-d/2}$. Getting back to the integral, we get
\[   \int_0^1\int_M\big|\operatorname{tr}(K_t(x,x))-p_t(x,x)\mathbb E^m_{t,x,x}[\mathbf1]\big|\mathrm dx \frac{\mathrm dt}{t^{1-s}}
\leq C\int_0^1\frac{\mathrm dt}{t^{1-s}}\cdot1\cdot t^2\cdot t^{-d/2}
   = C\int_0^1\frac{\mathrm dt}{t^{d/2-1-s}}. \]
This converges for all $s>d/2-2$.

\medskip

For $t$ large, we have
\[ \int_M\tr K_t(x,x)\mathrm dx-\frac1n\dim\big(\ker L\big)=\frac1n\sum_{\lambda\in\operatorname{sp}(L)\setminus\{0\}}\exp(-\lambda t)=O(\exp(-\delta t)), \]
for $\delta>0$ the smallest positive eigenvalue (because the eigenvalues grow polynomially). The integral in time is bounded as
\[   \int_1^\infty\frac{\mathrm dt}{t^{1-s}}\int_M\bigg|\operatorname{tr}K_t(x,x)-\frac{\dim\big(\ker L\big)}{n\cdot\operatorname{vol}(M)}\bigg|\mathrm dx
\leq C\int_1^\infty\frac{\exp(-\delta t)}{t^{1-s}}\mathrm dt, \]
which converges for all $s\in\mathbb R$.
\end{proof}

We can deduce the following: 
\begin{lemma}
    \label{le:zeta=E[chi]}
    Suppose $d<4$. Let $\nabla_0$ and $\nabla_1$ be two connections on two bundles $E_0$ and $E_1$ of rank $n_0$ and $n_1$, and assume that 
    \begin{equation} \label{eq:kernel}
    \frac{1}{n_0}\dim(\ker L_{\nabla_0})= 
    \frac{1}{n_1}\dim(\ker L_{\nabla_1}),
    \end{equation}
    which is automatically the case if $m$ is not identically vanishing or $M$ has a boundary, see Remark \ref{rk:kernel} below.
    
    Then,
    \[
    \frac{1}{n_1}\zeta_{\nabla_1}'(0)-\frac{1}{n_0}\zeta_{\nabla_0}'(0)
    = \int_0^\infty  \int_M \frac{p_t(x,x)}{t}  \mathbb{E}^m_{t,x,x}\big[\chi_{\nabla_0,\nabla_1}(W)\big] \d x \d t
    \]
    for $\chi_{\nabla_1,\nabla_0}$ the term arising from the difference in holonomy, see equation \eqref{eq:defXi}.
\end{lemma}

\begin{remark}
    The condition $d<4$ is necessary for the convergence of the integral close to zero. 
    Only in the specific case when $\nabla_0, \nabla_1$ are flat can the result be extended to higher dimensions. This specific case was treated in \cite{LeJan}.
    The same remark applies to the next section:
    if $d\geq4$ and $\nabla_0,\nabla_1$ are not flat, we do not even expect that their means converge unless in very specific cases.
\end{remark}

\begin{remark}
\label{rk:kernel}
To apply this result, it is useful to have a criterion ensuring that $L_\nabla$ has no kernel. Since this operator is elliptic, an element $s$ of the kernel must be smooth, and satisfy $L_\nabla s=0$ in the strong sense. Integrating by parts (the section must vanish on the boundary), we get
\[ \frac{1}{2}|\nabla s|_{L^2}^2 + |\sqrt m\,s|_{L^2}^2=\langle s,L_\nabla s\rangle=0. \]
This means that $\nabla s$ vanishes identically, i.e. that $s$ is parallel, and the value at any point is determined by its value at any other via the holonomy along any path.

If $m$ does not vanish identically or the boundary is non-empty, we can find a point where the section is zero, so by parallel transport it is zero everywhere. In fact, this condition depends on the mass but not on $\nabla$, and ensures that the kernel will be empty for \emph{all} connections on \emph{all} vector bundles.

In later sections, it will be crucial that not only $L_\nabla$ but also $\frac{1}{2} \Delta+m$ have a spectral gap, for this is what is needed for the Brownian loop soup to contain finitely many large loops. For this specific operator, and because $m$ is assumed to be pointwise non-negative, existence of a spectral gap is equivalent to either $\partial M$ being empty or $m$ being zero.
\end{remark}

\begin{remark}
\label{rk:kernel2}
The condition $M$ with boundary or $m$ non-vanishing are sufficient but not necessary for $L_\nabla$ to have a trivial kernel, and in fact there are actually very few cases where this is the case. Let us discuss two remaining obstructions.
 
The first is given by curvature issues. Let $s$ be an element of the kernel. For the purpose of this comment, the (global) \emph{holonomy group of $\nabla$ at $x$} is the closure of the set (easily seen to be a group) of all $\mathcal Hol^\nabla_{0,1}(\ell)\in O(E_x)$ for $\ell:[0,1]\to M$ a smooth loop based at $x$. Then $s(x)$ must be fixed by all the elements of this group, and a non-trivial kernel for $L_\nabla$ imposes some non-trivial conditions on the holonomy group. Let us define the infinitesimal, local and global holonomy groups $G_\mathrm{inf}\subseteq G_\mathrm{loc}\subseteq G_\mathrm{glob}$ at $x$ as the closed subgroups generated by the following elements. $G_\mathrm{inf}$ is generated by the elements of the form $\exp(R^\nabla(u,v))$, where $R^\nabla\in\mathcal C^\infty(M,\bigwedge^2T^*M\otimes\operatorname{End}(E))$ is the curvature form of $\nabla$ and $u,v$ range over the elements of the tangent space $T_xM$ at $x$; $G_\mathrm{loc}$ is generated by the holonomies of all smooth contractible loops based at $x$; $G_\mathrm{glob}$ is generated by the holonomies of all smooth loops based at $x$. Two generic elements of $O(E_x)$ will generate a dense subgroup when $n\geq2$ \cite[Lemma 8.1]{SurfaceGroups}, so generically $G_\mathrm{inf}$ is the full group when $d\geq3$ and $n\geq2$ or $d=n=2$; this is the case for instance if $E$ has rank two and the curvature is non-zero at $x$. To ensure that this does not happen, one can impose that the curvature be zero at $x$, which means that $G_\mathrm{inf}=\{\operatorname{Id}\}$. To go one step further, the Ambrose--Singer theorem \cite{AmbroseSinger} states that $G_\mathrm{loc}$ is described by the curvature at all points, in a very explicit manner. Without going into the details of this statement, let us simply say that it implies the following: the curvature must vanish at \emph{every} point of the manifold for $G_\mathrm{loc}$ to be trivial; in other words, the connection must be flat. Existence of a flat connection on a bundle is a non-trivial topological condition (for instance the bundle $T\mathbb S^2\to\mathbb S^2$ does not admit such a connection by the hairy ball theorem \cite[Chapter 5]{MilnorTopology}): $E$ must be a so-called flat bundle.

The second obstruction  arises from the global geometry of $M$. For $L_\nabla$ to have a trivial kernel, we need the long-range holonomies coming from the non-contractible loops to fix a line. This corresponds to some morphism $\pi_1(M)\to O(E_x)/G_\mathrm{loc}$ having a very small image (the image is $G_\mathrm{glob}/G_\mathrm{loc}$).

In conclusion, the existence of a non-trivial kernel forces $M$ to have no boundary, $m$ to vanish everywhere, the curvature of $\nabla$ to have a non-trivial kernel at every point, and the global holonomies to fix at least a line in this kernel, which is a rather strong collection of coincidences when the bundle has rank at least two.
\end{remark}

\begin{proof}[Proof of Lemma \ref{le:zeta=E[chi]}]
    We know that for $\Re z>d/2$,
    \begin{equation} 
    \label{eq:zetaIntegral}
    \frac1{n_1}\zeta_{\nabla_1}(z) - \frac1{n_0}\zeta_{\nabla_0}(z)
     = \frac1{\Gamma(z)}\int_0^\infty\frac{\mathrm dt}{t^{1-z}}\int_M\operatorname{tr}\Big(K^{\nabla_1}_t(x,x)-K^{\nabla_0}_t(x,x)\Big)\mathrm dx. \end{equation}
    The condition \eqref{eq:kernel} together with Corollary \ref{coro:integrability} allows to deduce that the right-hand side of \eqref{eq:zetaIntegral} is integrable whenever $\Re z>-1/2$. Besides, by monotonicity of both the integral over $[0,1]$ and the integral over $[1,\infty)$, this integrability actually holds locally uniform in $z$ in the half-plane $\Re z>-\frac{1}{2}$. Thus, the right-hand side of \eqref{eq:zetaIntegral} is holomorphic in this half-space. Since the left-hand side of \eqref{eq:zetaIntegral} is meromorphic in the same half-space, the equality extends to $\Re z>-\frac{1}{2}$, thus in particular around $z=0$. Noting that $(1/\Gamma)(0)=0$, $(1/\Gamma)'(0)=1$, we know that the left hand side in the lemma is equal to the integral of the right-hand side in \eqref{eq:zetaIntegral} with $z=0$.

    The result then follows from the representation of $K_t$ given by Theorem \ref{th:K=pE}.
\end{proof}

This concludes the first part of the proof. In the next section, we prove the second formula in Theorem \ref{th:main1}.

\section{Martingale convergence for some \texorpdfstring{multiplicative\\functionals}{multiplicative functionals} of a Brownian loop soup}
\label{sec:mart}
We now assume that either $m$ is not identically vanishing, or that $m=0$ but $M$ has a boundary. In both cases, the intensity measure $\Lambda_\delta$ on large loops is finite for all $\delta>0$; indeed, as discussed in Remark \ref{rk:kernel}, the first eigenvalue $\lambda_0$ of $\frac12\Delta+m$ must be positive, and 
\[ \int_M p_t(x,x) \mathbb{E}^m_{t,x,x}[\mathbf1]\d x=\Tr (e^{-t(\Delta/2+m)})=O(e^{-\lambda_0 t})\]
is integrable over $[\delta,\infty)$. It follows that $\mathcal{L}_\delta$ is a finite set for all $\delta>0$.

In this section we define an infinite product
\[ \prod_{\ell \in \mathcal{L}} \big(1+ \chi_{\nabla_0,\nabla_1}(\ell) \big), \]
for two connections $\nabla_0$, $\nabla_1$ over possibly different vector bundles, and fixed during this section. It must be understood that this product is improper, in that it is not expected to be absolutely convergent. We define it as the limit, as $\delta\to 0$, of the finite products 
\[ 
 Q_\delta \coloneqq \prod_{\ell \in \mathcal{L}_\delta} \big(1+ \chi_{\nabla_0,\nabla_1}(\ell) \big).
\]
One of our goal in this section is indeed to show that, for $d\leq 3$ and $m$ not identically zero, $Q_\delta$ is in $L^1$ for all $\delta>0$, and is convergent as $\delta\to 0$, both almost surely and in $L^1$.
We will do this by using the fact that the multiplicative increments of the process $(Q_\delta)_{\delta \searrow 0}$ are independent, from which it follows that the process $Q_\delta/ \mathbb{E}[Q_\delta]$ is a local martingale. From standard martingale theory, it then suffices to show that this local martingale is bounded in $L^2$ to deduce that it converges almost surely and in $L^2$. Proving that the expectation $\mathbb{E}[Q_\delta]$ also converges will conclude the proof. 

In the literature on the Brownian loop soup and its relation to determinant of Laplacians, additive rather than multiplicative functionals are considered. They require normalisation by substraction of logarithmically diverging counterterms, and it is then usually just the expectation of the sum that is shown to converge, rather than the sum itself.

Recall that we will prove the version of Theorem \ref{th:main1} only with a choice of intensity $\alpha=1$, the argument being, \emph{mutatis mutandis,} the same for other values.

\subsection{Variations around Campbell's theorem}

Our approach makes heavy use of some Campbell formula for multiplicative functionals of Poisson point processes. The version we give here is proved using the usual techniques for similar results, but we were not able to locate this precise statement in the literature. We isolate the proof in Appendix~\ref{app:Campbell}.

\begin{theorem}
\label{th:multiplicativeCampbell}
Let $S$ be a measurable space such that the diagonal is measurable in $S\times S$.
Let $\mathcal{P}$ be a Poisson process on $S$ with intensity measure $\mu$, and let $g:S\to \mathbb{C}$ be measurable.

If $g$ is in $L^1(S,\mu)$, then the product
\[ \Pi\coloneqq\prod_{x\in \mathcal P}\big(1+g(x)\big) \]
is absolutely convergent (in the sense that the sum of the $|g(x)|$ converges) almost surely, integrable, and
\[ \mathbb E[\Pi] = \exp\left(\int_Sg(x)\mu(\mathrm dx)\right). \]
\end{theorem}

Specifying the point process to our massive loop soup, we can rewrite the integrals appearing in the previous sections as expectations of infinite products.

\begin{lemma}
    \label{le:CampbellLoop}
    \hspace{-3.7pt}Let $\Phi:\mathscr{L}\to \mathbb{C}$ be an element of $L^0(\mathscr L,\Lambda)$ such that for all ${0<\delta<R<\infty}$,
    \[ \int_\delta^R\int_M  \mathbb{E}^m_{t,x,x}\big[|\Phi(W)|\big] \frac{p_t(x,x)}{t}\d x\d t <\infty. \]
    \begin{enumerate}
    \item \label{it:Campbellbulk}
        For almost all $(t,x)\in (0,+\infty)\times M$, $\Phi(W)$ is a well-defined random variable under $\mathbb{P}_{t,x,x}$.
        The function $(t,x)\mapsto \mathbb{E}_{t,x,x}[\Phi(W)]$ is measurable on $(0,+\infty) \times \mathbb{R}^2$, and for all $0<\delta<R<\infty$, 
        \begin{equation}
        \label{eq:integrabilityFunctionalbulk}
        \mathbb{E}\Big[ \prod_{\ell \in \mathcal{L}_\delta^R} (1+\Phi(\ell)) \Big]
        =\exp \Big( \int_\delta^R  \int_M \mathbb{E}^m_{t,x,x}[\Phi(W)]\frac{p_t(x,x)}{t} \d x\d t \Big).
        \end{equation}
    \item  \label{it:Campbellfunctional0}
        Assume that for some $R\in(0,\infty)$,
        \begin{equation} 
        \label{eq:integrabilityFunctional0}
        \int_0^R \int_M  |\mathbb{E}^m_{t,x,x}[\Phi(W)]| \frac{p_t(x,x)}{t}\d x\d t <\infty
        \end{equation}
       and the finite products $\prod_{\ell \in \mathcal{L}_\delta^R} (1+\Phi(\ell))$ converge in $L^1(\mathbb P)$ as $\delta\to 0$ to a limit that we denote by $\prod_{\ell \in \mathcal{L}^R} (1+\Phi(\ell)) $.
        Then, 
        \[ 
        \mathbb{E}\Big[ \prod_{\ell \in \mathcal{L}^R} (1+\Phi(\ell)) \Big]
        =\exp \Big( \int_0^R  \int_M \mathbb{E}^m_{t,x,x}[\Phi(W)] \frac{p_t(x,x)}{t} \d x\d t\Big).
        \]
    \item
    \label{it:CampbellfunctionalInfty}
        Assume that for some $\delta\in(0,\infty)$,
        \begin{equation}
        \label{eq:integrabilityFunctionalInfinity}
        \int_\delta^\infty \int_M \mathbb{E}^m_{t,x,x}[|\Phi(W)|] \frac{p_t(x,x)}{t}  \d x \d t<\infty.
        \end{equation}
        Then, the almost surely finite product $\prod_{\ell \in \mathcal{L}_\delta} (1+\Phi(\ell))$ is in $L^1(\mathbb P)$, and
        \[ 
        \mathbb{E}\Big[ \prod_{\ell \in \mathcal{L}_\delta} (1+\Phi(\ell)) \Big]
        =\exp \Big( \int_\delta^\infty \int_M \mathbb{E}^m_{t,x,x}[\Phi(W)] \frac{p_t(x,x)}{t} \d x \d t\Big).
        \]
    \end{enumerate}
\end{lemma}

\begin{remark}
\label{rem:positive}
For $\Phi:\mathscr L\to\mathbb C$ which satisfies the regularity and integrability assumptions of Lemma \ref{le:CampbellLoop} and such that additionally, $\Phi(\ell^{-1})=\overline{\Phi(\ell)}$ for $\ell^{-1}$ the time-reversal of $\ell$, it holds that for all $t,x$, $\mathbb{E}^m_{t,x,x}[\Phi(W)]$ is real (by invariance of $\mathbb P^m_{t,x,x}$ under time-reversal), and therefore the expectation \eqref{eq:integrabilityFunctionalbulk} is real-valued and positive.
\end{remark}

\begin{proof}[Proof of Lemma \ref{le:CampbellLoop}]
    Let $\Phi$ be as described. By our hypotheses on the $\sigma$-algebra of $\mathscr L$, we can see this map as defined on $(0,\infty)\times\mathcal C([0,1],M)$, and by Fubini's theorem $\Phi(W)$ is well-defined as a function of $L^0(\mathbb P_{t,x,x})$, for almost every $(t,x)$. Note also that on the product space $(0,\infty)\times\mathcal C([0,1],M)$, the topology is Hausdorff, so the diagonal is (closed hence) measurable in its square.

    We now apply the multiplicative version of Campbell's theorem to the Poisson point process $\mathcal{L}_\delta^R$, with the function $g=\Phi$.
    We obtain
   \begin{equation}
   \label{eq:temp1}
    \mathbb{E}\Big[ \prod_{\ell \in \mathcal{L}_\delta^R} (1+\Phi(\ell)) \Big]
    =\exp\Big( \int \Phi(\ell) \d \Lambda_\delta^R(\ell) \Big). 
    \end{equation}
    The disintegration formula 
    \[ \Lambda_\delta^R =\int_\delta^R \int_M \mathbb{P}^m_{t,x,x} \frac{p_t(x,x)}{t} \d x \d t \]
    directly implies that $\int \Phi(\ell) \d \Lambda_\delta^R(\ell)=\int_\delta^R \int_M \mathbb{E}^m_{t,x,x}[\Phi(W)] \frac{p_t(x,x)}{t} \d x \d t $, which  concludes the proof of the first point.

    The second point is obtained from the following equalities, justified below: 
    \begin{align*}
    \mathbb{E}\bigg[ \prod_{\ell\in \mathcal{L}^R} (1+\Phi(\ell))\bigg]
    &=\lim_{\delta\to 0}\mathbb{E}\bigg[ \prod_{\ell\in \mathcal{L}_\delta^R} (1+\Phi(\ell))\bigg]\\
    &= \lim_{\delta\to 0} \exp\Big( \int_\delta^R \int_M \mathbb{E}_{t,x,x}^m[\Phi(W)] \frac{p_t(x,x)}{t}\d t\d x\Big)\\
    &= \exp\Big( \int_0^R \int_M \mathbb{E}_{t,x,x}^m[\Phi(W)] \frac{p_t(x,x)}{t}\d t\d x\Big).
    \end{align*}
    The first equality follows from the assumption of convergence in $L^1(\Omega)$. The second equality follows from the first point of the lemma, that we have just proved. The third equality follows directly from the assumption \eqref{eq:integrabilityFunctional0}. This concludes the proof of the second point of the lemma.

    The third point follows again from the multiplicative version of Campbell's theorem: by the disintegration formula for $\Lambda_\delta$, the integrability assumption amounts to
    \[ 
    \int | \Phi(\ell)| \d \Lambda_\delta(\ell) <\infty.
    \]
    We can therefore apply the multiplicative Campbell's theorem directly to the Poisson point process $\mathcal{L}_\delta$, and to the function $g=\Phi$. We get
    \begin{align*}
    \mathbb{E}\Big[\prod_{\ell\in \mathcal{L} _\delta }(1+\Phi(\ell))\Big]
    &= \exp\Big( 
    \int \Phi(\ell) \d \Lambda_\delta (\ell)\Big)= \exp\Big( \int_\delta^\infty
    \int_M \mathbb{E}^m_{t,x,x}[\Phi(W)] \frac{p_t(x,x)}{t}\d x\d t\Big),
    \end{align*}
    where the last equality follows from the disintegration formula $\Lambda_\delta=\int_\delta^\infty \int_M \mathbb{P}^m_{t,x,x} \frac{p_t(x,x)}{t}\d x\d t $.
\end{proof}

\subsection{Convergence and expectation of the product over small loops}

From the last result, we can deduce the following:
\begin{lemma}
    \label{le:martingaleCV}
    Assume $d\leq 3$. 
    Recall that $\chi_{\nabla_0,\nabla_1}$ is defined by \eqref{eq:defXi}.  For any two smooth connections $\nabla_0,\nabla_1$ on bundles $E_0,E_1$ with possibly different ranks and for all $R\in(0,\infty)$, the quantity
    \[  Z_\delta^R\coloneqq \mathbb{E}\Big[ \prod_{\ell\in \mathcal{L}_\delta^R} \big(1+ \chi_{\nabla_0,\nabla_1}(\ell)\big)  \Big]>0
    \]
    converges as $\delta\to0$ toward a non-zero limit $Z^R$.
\end{lemma}

If we believe that the upper bound in the proof below is sharp in terms of powers of $t$, then it should be true that $Z^R_\delta$ diverges logarithmically when $d=4$, and as $\delta^{2-d/2}$ for $d>4$.

\begin{proof}
    The positivity of $Z^R_\delta$ follows from Remark \ref{rem:positive} with $\Phi=\chi$.
    Using \ref{it:Campbellbulk} of Lemma~\ref{le:CampbellLoop} ($\chi$ is bounded, so it satisfies the integrability property) then the technical Proposition~\ref{prop:tech} and the bound on $p$ given by \eqref{eq:unifpbound},
    we obtain, for any $0<{\delta'}<\delta<+\infty$,
    \begin{align}
    |\log(Z_{\delta'}^\delta)|
    &=\Big|\int_{\delta'}^\delta \int_M 
    \mathbb{E}_{t,x,x}^m[\chi(W)]
    \frac{p_t(x,x)}{t} \d t \d x\Big|\nonumber\\
    &\leq \int_{\delta'}^\delta \int_M 
    |\mathbb{E}^m_{t,x,x}[\chi(W)]|
    \frac{p_t(x,x)}{t} \d t \d x \nonumber\\
    &\leq C \int_{\delta'}^\delta\int_Mt^2\cdot\frac{t^{-d/2}}t\d x \d t, \label{eq:temp3} 
    \end{align}
    which converges as ${\delta'} \to 0$ since $d<4$.
    Thus, 
    \[ \lim_{\delta\to 0} \limsup_{{\delta'}\to 0} |\log(Z_{\delta'}^\delta)|=0.\]
    
    By the independence property of Poisson point processes, we deduce that for any $0<{\delta'}<\delta<R$,
    \[ \log(Z_{\delta'}^R)=\log(Z_{\delta'}^\delta)+\log(Z_\delta^R); \]
    by the above, for any sequence $(\delta_n)$ which decreases toward $0$, the sequence 
    $(\log(Z_{\delta_n}^R))$ must be a Cauchy sequence.
    If follows that $\log(Z _\delta^R)$ converges toward a finite limit as $\delta\to 0$, and therefore $Z _\delta^R$ converges toward a positive limit as $\delta\to 0$.
\end{proof}

\begin{lemma}
    \label{lem:convergenceofM}
    Assume $d\leq 3$. For $0<\delta\leq R<+\infty$, let \[M_\delta^R= \frac{1}{Z_\delta^R} \prod_{\ell\in \mathcal{L}_\delta^R} \big(1+ \chi_{\nabla_0,\nabla_1}(\ell)\big).\] 
    For all $R$, $\delta\mapsto M_\delta^R $ is a martingale as $\delta$ \emph{decreases,} with respect to the filtration $\delta\mapsto\sigma(\mathcal L_\delta)$. This martingale converges, almost surely and in $L^p(\mathbb{P})$ for all $p\geq 1$, as $\delta\to 0$, toward a limit $M^R$. 
\end{lemma}
\begin{proof}
    The fact that it is a martingale follows directly from the independence property of Poisson point processes and the fact that we have chosen the proper normalisation. 

    By Jensen inequality, it suffices to show the convergence holds in $L^{2k}$ for all positive integer $k$ in order to show it holds in $L^p$ for all $p$. Furthermore, by Doob's martingale convergence theorem, it actually suffices to show that it is bounded in $L^{2k}$ in order to show it converges in $L^{2k}$ and almost surely. Since the constants $Z^R_\delta$ converges toward a non zero value as $\delta\to 0$, it actually suffices to show that $    Q_\delta^R \coloneqq Z^R_\delta M^R_\delta$ is bounded in $L^{2k}$ for all positive integer $k$. 

    Since    
    \[
    |Q_\delta^R|^{2k} = \prod_{\ell\in \mathcal{L}_\delta^R} |1+ \chi(\ell)|^{2k},
    \]
    we use \ref{it:Campbellbulk} of Lemma \ref{le:CampbellLoop} again, this time with 
    \[\Phi:W\mapsto \big|1+ \chi(W)\big|^{2k}-1, \]
    to get
    \begin{align*}
    \mathbb{E} [ |Q_\delta^R|^{2k}]
    &=\exp\Big( \int_\delta^R \int_M  \mathbb{E}^m_{t,x,x}[|1+\chi(W)|^{2k}-1]    \frac{p_t(x,x)}{t}\d x \d t \Big)\\
    &=\exp\Big( \int_\delta^R \int_M  \sum_{j=0}^{k-1} \mathbb{E}^m_{t,x,x}[(|1+\chi(W)|^{2}-1)|1+\chi(W)|^{2j}  ]    \frac{p_t(x,x)}{t}\d x \d t \Big)\\
    &\leq \exp\Big( \int_\delta^R \int_M  \sum_{j=0}^{k-1} \mathbb{E}^m_{t,x,x}[|\chi(W)|^{2}+2|\mathfrak{R}(\chi(W)) | )3^{2j}  ]    \frac{p_t(x,x)}{t}\d x \d t \Big)\\
    &\leq \exp\Big( \int_\delta^R \int_M  3^k C t^2     \frac{p_t(x,x)}{t}\d x \d t \Big).
    \end{align*}
   The last inequality is obtained as follows. For a unitary $U$, $\mathfrak{R}(\tr(U))=\frac{1}{2}\mathfrak{R}( \tr(U+U^*))$. Hence, for $U_0,U_1$ unitary (on possibly different vector spaces), 
    \begin{align*}
    |\mathfrak{R}(\tr(U_1)-\tr(U_0))|
    &= |\mathfrak{R}(\frac{1}{2}\tr(2 \operatorname{Id}- U_1-U_1^*)-\frac{1}{2}\tr(2 \operatorname{Id}- U_0-U_0^*))|\\
    &\leq \frac{1}{2}|2I-U_1-U_1^* |+|2I-U_0-U_0^*|=
    \frac{1}{2}(|I-U_1|^2 +|I-U_0|^2).
    \end{align*}
    Using this equality with $U_i=\mathcal Hol^{\nabla_i}(W)$ and using the technical proposition \ref{prop:tech},  
    we deduce     
    \[    \mathbb{E}^m_{t,x,x}[|\chi(W)|^{2}+2|\mathfrak{R}(\chi(W)) |]\leq 2 Ct^2 \]
    for some $C$.
    From \eqref{eq:unifpbound}, the kernel is of order $t^{-d/2}$ as $t\to 0$, and it follows that for some constant $C'>0$,
    \begin{align*}
    \mathbb{E} [ |Q_\delta^R|^{2k}] \leq \exp
    \Big( C' 3^k \int_0^R \int_M t^2\cdot\frac{t^{-d/2}}t \d x \d t \Big)<+\infty.
    \end{align*}
    Notice that the finiteness of this last integral follows again from $d\leq 3$. This concludes the proof.
\end{proof}

\begin{corollary}
    \label{coro:smallloops}
    For all $R<+\infty$, the products $Q_\delta^R$ converge, almost surely and in $L^p(\mathbb{P})$ for all $p\geq1$, as $\delta\to 0$, toward 
    $ Q^R= Z^R M^R$.

    Furthermore,
    \[ \mathbb{E}[ Q^R]=\exp\Big( 
    \int_0^R\int_M \mathbb{E}^m_{t,x,x}[ \chi_{\nabla_0,\nabla_1}(W) ] \frac{p_t(x,x)}{t} \d x\d t
    \Big).  \]
\end{corollary}
\begin{proof}
    The almost sure (resp. $L^p$) convergence of $Q_\delta^R$ follows from the convergence of $Z_\delta^R$ and the almost sure (resp. $L^p$) convergence of $M_\delta^R$: 
    \[ \|Q_\delta^R-Q^R\|_{L^p}=\|Z_\delta^RM_\delta^R-Z^RM^R\|_{L^p}\leq Z_\delta^R\| M_\delta^R-M^R\|_{L^p}+|Z_\delta^R-Z^R|\|M^R\|_{L^p}\underset{\delta\to 0}\longrightarrow 0. \]
    
    To prove the second point, we apply \ref{it:Campbellfunctional0} of Lemma \ref{le:CampbellLoop} with $\Phi=\chi$. 
    We have already checked that
    \[ \int_0^R \int_M \big|\mathbb{E}^m_{t,x,x}[\Phi(W)]\big|\frac{p_t(x,x)}{t} \d x\d t <\infty    \]
    (see Equation \eqref{eq:temp3}), and we have just checked that the product $Q_\delta^R=\prod_{\ell \in \mathcal{L}_\delta^R} (1+\Phi(\ell))$ converges in $L^1(\Omega)$. The conclusion of \ref{it:Campbellfunctional0} of Lemma \ref{le:CampbellLoop} gives the desired expression for $\mathbb{E}[ Q^R]$. 
\end{proof}

\subsection{Adding the large loops}
\begin{lemma}
    \label{le:bigloops}
    Assume $m\neq 0$ or $\partial M\neq \emptyset$. Then, for all $\delta\in(0,\infty)$ and $p\geq 1$, the finite product 
    $Q_\delta \coloneqq \prod_{W \in \mathcal{L}_\delta} (1+\chi_{\nabla_0,\nabla_1}(W) ) $ is in $L^p(\mathbb{P})$, and 
    \[ \mathbb{E}[Q_\delta]=\exp\Big( \int_\delta^\infty \int_M  \mathbb{E}^m_{t,x,x}[\chi_{\nabla_0,\nabla_1}(W)]\frac{p_t(x,x)}{t}  \d x\d t\Big). \]
\end{lemma}
\begin{proof}
    To show $\mathbb{E}[|Q_\delta|^{p}]<\infty$, we apply \ref{it:CampbellfunctionalInfty} of Lemma \ref{le:CampbellLoop} with $\Phi(W)= |1+\chi(W)|^{p}-1$. Since $|1+\chi(W)|^{p}-1\leq 3^p-1$, 
    \[\int_\delta^\infty \int_M  \mathbb{E}^m_{t,x,x}[  |\Phi(W)| ] \frac{p_t(x,x)}{t}  \d x\d t\leq  (3^p-1) \int_\delta^\infty \Tr( e^{-t (\Delta/2+m) } ) \d t<\infty, 
    \]
    as we have seen at the beginning of this section.

    The expression for $\mathbb{E}[Q_\delta]$ then follows from \ref{it:CampbellfunctionalInfty} of Lemma \ref{le:CampbellLoop}  with $\Phi(W)= \chi(W)$.
\end{proof}

We immediately deduce the following, which concludes the proof of Theorem \ref{th:main1} (conditional on the technical estimation, Proposition \ref{prop:tech}).
\begin{prop}
\label{prop:loops}
  Assume $d\leq 3$ and either $m\neq 0$ or $\partial M\neq \emptyset$. Then, for all $p\geq 1$, $Q_\delta\in L^p(\mathbb{P})$ for all $\delta>0$, and $Q_\delta$ converges almost surely and in $L^p$ as $\delta\to 0$. The limit $Q$ satisfies 
  \[ \mathbb{E}[Q]=\exp\Big( \int_0^\infty \int_M \mathbb{E}^m_{t,x,x} [\chi_{\nabla_0,\nabla_1}(W)] \frac{p_t(x,x)}{t}\d x \d t\Big). \]
\end{prop}
\begin{proof}
    The fact $Q_\delta\in L^p$ follows from Lemma \ref{le:bigloops}. 
    Since $Q_\delta=Q_\delta^1 Q_1$ for all $\delta\leq 1$ and since $Q_\delta^1$ converges toward $Q^1$, almost surely and in $L^{2p}$, (by Corollary \ref{coro:smallloops}, and since $Q_1\in L^{2p}$, we deduce $Q_\delta$ converges almost surely and in $L^p$ toward $Q\coloneqq Q_1Q^1$ as $\delta\to 0$. All that remains to be shown is the expression for $\mathbb{E}[Q]$. Since $Q_1$ is $\mathcal{L}_1$-measurable whilst $Q^1$ is $\mathcal{L}^1$-measurable, $Q_1$ and $Q^1$ are independent. Using the expressions for $\mathbb{E}[Q_1]$ and $\mathbb{E}[Q^1]$ given by  Corollary \ref{coro:smallloops} and Lemma \ref{le:bigloops}, we get 
    \begin{align*} 
    \mathbb{E}[Q]&=\mathbb{E}[Q_1]\mathbb{E}[Q^1]\\
    &= \exp\Big( \int_0^1 \int_M \mathbb{E}^m_{t,x,x}[\chi(W)] \frac{p_t(x,x)}{t} \d x \d t\Big)\exp\Big( \int_1^\infty \int_M \mathbb{E}_{t,x,x}^m[\chi(W)] \frac{p_t(x,x)}{t} \d x \d t\Big)\\
    &=\exp\Big( \int_0^\infty \int_M \mathbb{E}^m_{t,x,x} [\chi(W)] \frac{p_t(x,x)}{t}\d x \d t\Big),
    \end{align*}
    which concludes the proof.
\end{proof}

\begin{remark}
    As we remarked in Corollary \ref{cor:main}, we can compare any connection $\nabla$ on any bundle $E$ to the flat connection $\mathrm d$ on the trivial bundle $M\times\mathbb R$, and from equation \eqref{eq:pre-ratio} of Theorem \ref{th:main1} we immediately deduce
    \[    \exp\left( -\frac{\zeta'_\nabla(0)}{n} +\zeta'_{\d}(0)  \right) = \mathbb E\bigg[\prod_{\ell\in\mathcal L}\tr(\mathcal Hol_\nabla(\ell))\bigg]^{-1} \geq 1.\]
    The inequality
    \begin{equation}
    \label{eq:diamagnetic}
    \Det(L_\nabla^{-1}) = \exp\left(\zeta_\nabla'(0)\right) \leq \exp(\zeta_{\mathrm d}'(0))^n,
    \end{equation}
    which we will use in section \ref{ssec:minimalcouplingGFF} below, was first proved in \cite{Schrader} in the massless case (with a boundary) using geometric analysis, but is a free by-product of our probabilistic representation. Note that this bound depends on the bundle only through its rank.
\end{remark}

\section{Proof of the technical estimate}
\label{sec:tech}

In this section, we prove the estimate in Proposition \ref{prop:tech}, i.e. we derive an explicit $L^p$ bound on the holonomy along Brownian loops, which we used above for small times. We proceed in three steps: we show first that the first estimate of the proposition holds for $M$ a $d$-dimensional ball (with a generic metric), then we extend this bound for all manifolds, and finally we deduce the three other estimates of Proposition \ref{prop:tech}. We work in the slightly more general setting of principal bundles with compact fiber. Since a Euclidean (resp. Hermitian) metric bundle $E$ is the vector bundle associated to a principal bundle with fiber $O_n(\mathbb R)$ (resp. $U_n(\mathbb C)$), namely its frame bundle, we will be able to reduce from it the result we want.

Recall that we work on a connected manifold $M$ of dimension $d$, endowed with a metric $g$. Fixing some compact Lie group $G$, we work with a principal $G$-bundle $P$ over $M$. Writing $\mathfrak g$ for the Lie algebra of $G$, we define the adjoint bundle $P_{\mathfrak g}$ of $P$, whose underlying set is the quotient of $P\times\mathfrak g$ under the action
\[ g\cdot(p,v)\coloneqq(p\cdot g^{-1},gvg^{-1}) \]
of $G$. Note that in general $P_{\mathfrak g}\not\simeq M\times\mathfrak g$ (a phenomenon similar to $\operatorname{End}(E)\not\simeq M\times\operatorname{End}(\mathbb R^n)$ for $E$ a non-trivial vector bundle of rank $n$).

Without loss of generality, we assume that $G$ is a closed subgroup of $O_n(\mathbb R)$ for some $n\geq1$. We also use Einstein notation (summing over repeated indices) throughout.
\medskip

We conclude this small introduction by a technical estimate about general Stratonovich integrals. Before we introduce it, let us describe in informal terms the context in which we aim to apply it.

Let $W$ be a Brownian motion in $M$ with respect to the metric $g$, and $Z$ a process over $W$, in the sense that it takes values in a bundle over $M$ and its projection is $W$; think for instance of the holonomy along $W$, or of the integral of the mass $m$. The processes we have in mind are functionals of the type $t\mapsto\alpha(Z_t)\in\mathbb R$ (think one coordinate of $Z$, in a local system). Namely, in the following lemma, a direct application of It\^o's formula shows that we can choose $X$ and $Y$ of such a form $t\mapsto\alpha(Z_t)$, provided that $Z$ be the solution to a stochastic differential equation with smooth coefficients driven by $W$, $\alpha$ be smooth, and $\tau$ be a stopping time such that $Z_{|[0,\tau)}$ takes values in a compact subset of the bundle.
\begin{lemma}
\label{lem:StratoLp}
Let $X$ and $Y$ be two one-dimensional continuous semimartingales defined at least up to some stopping time $\tau$. Suppose that the brackets of $X$ and $Y$ admit finite Lipschitz constants $\mu_X^2$ and $\mu_Y^2$, and suppose also that the finite-variation part of $Y$ admits the finite Lipschitz constant $\nu_Y$. There exists a universal constant $C$ such that, for all $t\in[0,1]$ and $p\geq2$,
\[ \mathbb E\left[\left|\int_0^{t\wedge\tau}X_s\circ\mathrm dY_s\right|^p\right]\leq \big(C\mu_X\mu_Y\big)^pt^p + \big(C(\mu_Y+\nu_Y)\big)^pp^{p/2}t^{p/2-1}\int_0^t\mathbb E\big[|X_s|^p\mathbf1_{s<\tau}\big]\mathrm ds. \]
In particular, in the case $X\equiv1$ and $Y_0=0$,
\[ \mathbb E\big[|Y_{t\wedge\tau}|^p\big] \leq \big(C(\mu_Y+\nu_Y)\big)^pp^{p/2}t^{p/2}. \]
\end{lemma}

\begin{proof}
Let us write $A$ and $M$ for the finite-variation and local martingale parts of $Y-Y_0$. We know that
\[ \int_0^{t\wedge\tau}X_s\circ\mathrm dY_s
 = \int_0^{t\wedge\tau}X_s\mathrm dA_s
 + \int_0^{t\wedge\tau}X_s\mathrm dM_s
 + \frac12\langle X,Y\rangle_{t\wedge\tau}
=: a_t + b_t + c_t, \]
so by the Hölder inequality $(a+b+c)^p\leq 3^{p-1} (a^p+b^p+c^p)$,
\[   \mathbb E\Big[\Big|\int_0^{t\wedge\tau}X_s\circ\mathrm dY_s\Big|^p\Big]
\leq 3^{p-1}\big(\mathbb E\big[|a_t|^p\big] + \mathbb E\big[|b_t|^p\big] + \mathbb E\big[|c_t|^p\big]\big). \]

We know that $c_t^2\leq\langle X,X\rangle_{t\wedge\tau}\langle Y,Y\rangle_{t\wedge\tau}$ (because the quadratic form $(\alpha,\beta)\mapsto\langle\alpha X+\beta Y,\alpha X+\beta Y\rangle_{t\wedge\tau}$ is semi-definite positive), so
\[   \mathbb E\big[|c_t|^p\big]
\leq \mathbb E\Big[\big((\mu_X^2t)^{1/2}(\mu_Y^2t)^{1/2}\big)^p\Big]=(\mu_X\mu_Yt)^p. \]
For the finite variation part, we get
\[   \mathbb E\big[|a_t|^p\big]
\leq \mathbb E\left[\left(\int_0^{t\wedge\tau}|X_s|\nu_Y\mathrm ds\right)^p\right]
\leq t^{p-1}\nu_Y^p\int_0^t\mathbb E\big[|X_s|^p\mathbf1_{s<\tau}\big]\mathrm ds. \]
For the martingale part, we use the Burkholder--Davis--Gundy inequality for continuous martingales:%
\footnote{
The usual constant found in references takes the form $(Cp)^p$. 
The optimal constant for \emph{continuous} martingales was found by Davis \cite[Theorem 3.1]{Davis} in terms of the largest $0$ of some Hermite polynomial, and this value is known to be smaller than  $C^p p^{p/2}$ for an explicit $C$  (see \cite[Equation (6.2.18)]{Szego} for a more precise statement). See also \cite{zakai} for a more direct result that is far sufficient for our purpose. 
}
\begin{align*}
      \mathbb E\big[|b_t|^p\big]
&\leq C^pp^{p/2}\mathbb E\left[\left(\int_0^{t\wedge\tau}|X_s|^2\mathrm d\langle Y,Y\rangle_s\right)^{p/2}\right] \\
&\leq (C\mu_Y)^pp^{p/2}t^{p/2-1}\int_0^t\mathbb E\big[|X_s|^p\mathbf1_{s<\tau}\big]\mathrm ds.
\end{align*}

The first inequality of the lemma follows putting together those three bounds, noting that $t^p\leq p^{p/2}t^{p/2}$. The second is immediate, since we can choose $\mu_X=0$.
\end{proof}

\subsection{Local estimates}

In this section, we assume that the manifold $M$ is the ball $B_0(2)\subset\mathbb R^d$, which later will be thought of as a neighbourhood of a point in a general manifold. The metric $g$ we consider, however, is \emph{not} the Euclidean metric; it is an arbitrary smooth Riemannian structure. We will say ``unit ball'' to refer to the Euclidean ball $B_0(1)$, together with the restriction of $g$; it is just a question of convenience, since we could just as well use any relatively compact open neighbourhood of 0.

Since our base space is contractible, the topology plays no role here, and all bundles are trivial. We assume, up to isomorphism, that $P=M\times G$ and $P_{\mathfrak g}=M\times\mathfrak g$.
It will be easier to work with vector quantities, so we define $E$ the trivial bundle $M\times\mathbb R^n$. It is associated to $P$, and a connection $\omega$ on $P$ gives rise to a connection $\nabla=\nabla^\omega$ on $E$. For $\gamma:[0,t]\to M$ a regular curve, the parallel transport $\mathcal Hol^\nabla_{0,t}(\gamma):E_{\gamma_0}\to E_{\gamma_t}$ in $E$ and the holonomy $\mathcal Hol^\omega_{0,t}(\gamma)\in G$ (the only curve starting from the identity such that $s\mapsto(\gamma_s,\mathcal Hol^\omega_{0,s}(\gamma)\cdot h)$ is parallel in $P$ for all $h\in G$) are related through
\[ \mathcal Hol^\nabla_{0,t}(\gamma)(v)=\mathcal Hol^\omega_{0,t}(\gamma)\cdot v. \]
Note that the holonomy of $\omega$ as an element of $G$ makes sense only because $P$ is an explicit product; otherwise, it would be a map from $P_{\gamma_0}$ to $P_{\gamma_t}$ such that $\mathcal Hol^\omega_{0,t}(\gamma)(p\cdot g)=\mathcal Hol^\omega_{0,t}(\gamma)(p)\cdot g$ for all $g\in G$, $p\in P_{\gamma_0}$; it is also determined by the image of a single element, but there is no identity element in $P_{\gamma_0}$ to play the role of a basepoint. It is well-defined again if $P$ is not trivial but $\gamma$ is a loop.

\medskip

Our first step is to construct convenient flat connections; indeed, as discussed in the remark after the statement of Proposition \ref{prop:tech}, a flat connection does not contribute (locally) to the holonomy of a \emph{loop,} so we can add or subtract it to the holonomy, but it can cancel the main contribution of the holonomy of an \emph{open path.} This is the point of the condition $A_x(x)=0$ in the coming lemma.

Let $D$ be the connection on all our objects coming from the product structure. If needed, we can specify the bundle as a superscript, e.g. $D^{TM}$ or $D^E$. Note that the difference between two connections on a vector bundle $V$ is generally identified with a section of the bundle $T^*M\otimes\mathrm{End}(V)$, but if $V$ is associated to $P$ and both of these connections actually come from connections on $P$, then the corresponding endomorphisms of $V$ come from the action of the Lie algebra $\mathfrak g$, and the difference lies within the sections of the bundle $T^*M\otimes P_{\mathfrak g}$. For everything we consider in the following, we will then be able to write principal connections $\nabla$ over $E$ as a sum $D^E+(B\cdot{})$, for some $B$ with values in $T^*M\otimes P_{\mathfrak g}$, and a perturbation of a connection $\nabla$ will also be of the form $\nabla+(A\cdot{})$.

\begin{lemma}
\label{lem:flatmodification}
Let $\nabla$ be a connection on $E$ coming from a (smooth) connection on $P$. There exists a smooth section $A:(x,y)\mapsto A_x(y)\in T^*_yM\otimes\mathfrak g$ of the vector bundle
\begin{align*}
M\times\big(T^*M\otimes P_{\mathfrak g}\big)&\to M\times M \\
(x,\alpha)&\mapsto\big(x,\pi(\alpha)\big)
\end{align*}
over $M\times M$ such that $\nabla+(A_x\cdot{})$ is flat and $A_x(x)=0$ for all fixed $x\in M$.
\end{lemma}

\begin{proof}
As discussed above, $\nabla-D$ can be identified with a one-form $B$ over $M$ with values in $\mathfrak g$. Define the section
\[ \Phi:(x,y)\mapsto\exp\big(-B(x)(y-x)\big) \]
of the bundle $M\times P\to M\times M$, and define $A$ through the relation
\[ A_x(y)
 = -B(y) - \mathrm dx^i\otimes\left(\frac{\mathrm d}{\mathrm dy_i}\Phi_x(y)\right)\Phi_x(y)^{-1}\in T^*_yM\otimes\mathfrak g\simeq(T^*M\otimes P_{\mathfrak g})_y. \]
Clearly $A$ is smooth in $(x,y)$, and a simple computation yields $A_x(x)=0$. We will show that for all $x$ the connection $\nabla+(A_x\cdot{})$ is flat by exhibiting horizontal sections that induce a basis of $E_y$ at every $y\in M$.

Fix $x\in M$. For all $v\in\mathbb R^n$, define the section
\[ s_v:y\mapsto \Phi_x(y)\cdot v \]
of $E$. Clearly $s_v(x)=v$, so the $(s_v(x))_{v\in\mathbb R^n}$ generate $E_x$. Moreover,
\begin{align*}
   \big(\nabla s_v + A_x\cdot s_v\big)(y)
&= \mathrm dx^i\otimes\Big((\partial_is_v)(y)
 + \big(B(y)_i\cdot s_v(y)\big)
 + \big(A_x(y)_i\cdot s_v(y)\big)\Big) \\
&= \mathrm dx^i\otimes\left(\frac{\mathrm d}{\mathrm dy_i}\Phi_x(y)\cdot v\right) \\
&\qquad - \mathrm dx^i\otimes\left(\frac{\mathrm d}{\mathrm dy_i}\Phi_x(y)\right)\Phi_x(y)^{-1}\cdot s_v(y) \\
&= 0,
\end{align*}
so $s_v$ is parallel with respect to $\nabla+(A_x\cdot{})$, for our choice of fixed $x$. In particular, $\nabla+(A_x\cdot{})$ is indeed flat for every fixed $x\in M$.
\end{proof}

The notion of parallel transport, seen as a linear map from the corresponding fibers of $E$, that we used so far for smooth curves and Brownian paths, is well-defined (in the Stratonovich sense again) over any semimartingale by the standard stochastic calculus. We write $\mathbb P_x$ for the probability measure, on the space of continuous curves on $M$, for which the canonical variable $W$ (the identity) becomes a Brownian motion started from $x$, with respect to the metric $g$ on $M=B_0(2)$.

The next result is precisely the local form of the estimate we seek.

\begin{lemma}
\label{lem:topologylessholonomyestimate}
Let $\nabla$ and $A$ be as in Lemma \ref{lem:flatmodification}, or more generally $A$ is such that $A_x(x)$ vanishes for all $x\in M$ but $\nabla+(A\cdot{})$ may not satisfy the flatness condition. Let $\tau$ be the exit time of the unit ball for $W$.
There exists a constant $C=C_{\nabla,A,g}>0$ such that
\[   \sup_{\substack{\nabla'=\nabla+(A_x\cdot{})\\x\in B_0(1)}}\mathbb E_x\big[|\mathcal Hol^{\nabla'}_{0,t\wedge\tau}(W)-\mathcal Hol^\nabla_{0,t\wedge\tau}(W)|^p\big]
\leq (Ctp)^p \]
for all $t\geq0$ and $p\geq1$.
\end{lemma}

The norm we choose for the operators does not matter since the constants are not made explicit, but for the sake of convenience we take it to be the operator norm. The various other tensors can be endowed with any norm, as far as this proof is concerned.

\begin{proof}
We fix $x\in B_0(1)$ and set $\nabla'=\nabla+(A_x\cdot{})$. Using Hölder's inequality, we reduce to the case $p\geq2$. We write $C$ for a constant independent of $x$, $t$ and $p$, allowed to vary from one inequality to the other. We first consider the case $t\in[0,1]$.

The difference $\nabla-D$ can be identified with a one form $B$ over $M$ with values in $\mathfrak g$. By definition of the parallel transport, we must have for all $v\in E_x$ and $t<\tau$ (i.e. the corresponding integral equation holds for all $t$ when the upper bound is $t\wedge\tau$)
\[ 0
 = \mathrm d^\nabla\big(\mathcal Hol^\nabla_{0,t}(W)v\big)
 = \Big(\mathrm d\big(\mathcal Hol^\nabla_{0,t}(W)_\beta^\alpha\big)\Big)v^\beta\epsilon_\alpha
 + B(W_t)\big({}\circ\mathrm dW_t\big)\mathcal Hol^\nabla_{0,t}(W)v. \]
In other words,
\begin{equation}\label{eq:SDEHol}
   \mathrm d(\mathcal Hol^\nabla_{0,t}(W)_\beta^\alpha\big)
 = -\big(B(W_t)\big)_{i\nu}^\alpha\mathcal Hol^\nabla_{0,t}(W)_\beta^\nu\circ\mathrm dW_t^i, \end{equation}
and a similar relation holds for the parallel transport with respect to $\nabla'$, where among others $B$ must be replaced by $B+A_x$.

Setting
\[ Y:t\mapsto\mathcal Hol^{\nabla'}_{0,t\wedge\tau}(W)-\mathcal Hol^\nabla_{0,t\wedge\tau}(W), \]
we get, omitting indices,
\begin{align*}
   Y_t
&= - \int_0^{t\wedge\tau}A_x(W_s)\big({}\circ\mathrm dW_s\big)\mathcal Hol^{\nabla'}_{0,s}(W)
   - \int_0^{t\wedge\tau}B(W_s)\big({}\circ\mathrm dW_s\big)Y_s \\
&= -\int_0^{t\wedge\tau}\Phi_s\circ\mathrm dW_s - \int_0^{t\wedge\tau}\big({}\circ\mathrm d\Psi_s)Y_s,
\end{align*}
where
\begin{align*}
(\Phi_t)_{i\alpha}^\beta&=\big(A_x(W_{t\wedge\tau})\big)_{i\nu}^\beta\big(\mathcal Hol^{\nabla'}_{0,t\wedge\tau}(W)\big)^\nu_\alpha, &
(\Psi_t)_\alpha^\beta=\int_0^{t\wedge\tau}\big(B(W_s)\big)_{i\alpha}^\beta\circ\mathrm dW_s^i.
\end{align*}

It should be clear that all the coefficients of $\Phi$, $W$, $\Psi$ and $Y$ satisfy the hypotheses of the stochastic analysis lemma \ref{lem:StratoLp}, in the sense that their brackets and finite variation parts are uniformly Lipschitz; indeed, they are of the form described above the statement of the lemma, where the process over $W$ is
\[ t\mapsto\big(W_t,\Psi_t,\mathcal Hol^{\nabla'}_{0,t}(W),\mathcal Hol^\nabla_{0,t}(W)\big). \]
Accordingly, we find for $t\in[0,1]$, which we now assume,
\[   \mathbb E\big[|Y_t|^p\big]
\leq \big(Ct)^p
   + C^pp^{p/2}t^{p/2-1}\int_0^t\Big(\mathbb E\big[|\Phi_s|^p\mathbf1_{s<\tau}\big]+\mathbb E\big[|Y_s|^p\mathbf1_{s<\tau}\big]\Big)\mathrm ds. \]
We also have, using the second inequality of the lemma (note that $\Phi_0=0$),
\[   \int_0^t\mathbb E\big[|\Phi_s|^p\mathbf1_{s<\tau}\big]\mathrm ds
\leq \int_0^t\mathbb E\big[|\Phi_{s\wedge\tau}|^p\big]\mathrm ds
\leq \int_0^tC^pp^{p/2}s^{p/2}\mathrm ds
\leq C^pp^{p/2}t^{p/2+1}. \]
This gives the bound on $Y$
\[   \mathbb E\big[|Y_t|^p\big]
\leq \big(Cpt)^p
   + C^pp^{p/2}t^{p/2-1}\int_0^t\mathbb E\big[|Y_s|^p\big]\mathrm ds, \]
which gives the following control by Grönwall's lemma:
\[   \mathbb E\big[|Y_t|^p\big]
\leq (Cpt)^p\exp\big((Cpt)^{p/2}\big). \]
Recall that this inequality is valid for all $x\in B_0(1)$, $t\in [0,1]$, $p\geq 2$. Without loss of generality, we choose $C$ larger than $1$.

We now consider a general $t\geq0$. Since the holonomies are unitary, we have $\mathbb E\big[|Y_t|^p\big]\leq2^p $, and the result follows from the fact that
\[ 2^p\wedge\Big((Cpt)^p\exp((Cpt)^{p/2})+\infty\mathbf1_{t>1}\Big)\leq \big(2Cpt\big)^p \]
over our expected range of $p$ and $t$, which one can show by arguing whether $Cpt\geq1$ or not.
\end{proof}

Finally, to compare the local model to the global situation, we need to estimate the probability that we leave a neighbourhood of the starting point. This bound is probably known to the experts but we were not able to locate it in the literature, so we present a short proof based on a quantitative estimate we learned from Hsu.

\begin{lemma}
\label{lem:exitball}
Let $\tau$ be the exit time of the Euclidean unit ball $B_0(1)\subset\mathbb R^d$ for the Brownian motion associated to the metric $g$. There exists a constant $C=C_g>0$ such that
\[ \sup_{x\in B_0(1/2)}\mathbb P_x(\tau\leq t)\leq C\exp(-1/(Ct)) \]
for all $t>0$.
\end{lemma}

\begin{proof}
Choose some $\delta>0$ that is smaller than the $d_g$-distance between $\overline B_0(1/2)$ and $B_0(1)^\complement$. Under $\mathbb P_x$, the process
\[ \widetilde W:t\mapsto W_{\delta^2t} \]
is a Brownian motion associated to $(M,g/\delta^2)$. The Ricci curvature of $(M,g/\delta^2)$ is the same as that of $(M,g)$. Let $-K^2<0$ be a lower bound on the Ricci curvature of $g$ over all of $B_0(1)$, so that it is also a lower bound on the Ricci curvature of $(B_0(1),g/\delta^2)$. Without loss of generality, we choose $K\geq1$.

For all $x\in B_0(1/2)$, $t>0$,
\begin{align*}
      \mathbb P_x(\tau<t)
&\leq \mathbb P_x\left(\sup_{s\leq t}d_g(W_0,W_{t\wedge\tau})>\delta\right)
   = \mathbb P_x\left(\sup_{s\leq t/\delta^2}d_{g/\delta^2}\big(\widetilde W_0,\widetilde W_{t\wedge(\tau/\delta^2)}\big)>1\right).
\end{align*}
According to \cite[Theorem 3.6.1]{Hsu}, there exists some $\eta>0$ depending only on $d$ such that whenever $L\geq1$ and $-L^2$ is a bound on the Ricci curvature,
\[   \mathbb P_x\left(\sup_{s\leq\eta/L}d_{g/\delta^2}\big(\widetilde W_0,\widetilde W_{t\wedge(\tau/\delta^2)}\big)>1\right)
\leq \exp(-L/2). \]
In particular, for all $t\leq\delta^2\eta/K$, we can choose $L=\delta^2\eta/t$ and get
\[   \mathbb P_x(\tau<t)
\leq \mathbb P_x\left(\sup_{s\leq t/\delta^2}d_{g/\delta^2}\big(\widetilde W_0,\widetilde W_{t\wedge(\tau/\delta^2)}\big)>1\right)
\leq \exp\big(-\delta^2\eta/2t\big). \]
Up to adding a large multiplicative factor to the exponential term, and since $\eta$, $\delta$ and $K$ are independent of $x$, we get the expected inequality for all $t>0$.
\end{proof}

\subsection{Proof of Proposition \ref{prop:tech}}

We now go back to the general setting of the section, i.e. $M$ is a closed manifold rather than a topological ball. Recall that $g$ is a metric on $M$, $G$ a closed subgroup of $O_n(\mathbb R)$, $P$ a principal $G$-bundle over $M$. We consider the holonomy of some principal connection $\omega$ on $P$, from which all our vector connections $\nabla$ derive.

\begin{theorem}
\label{th:momentsoftheholonomy}
There exists some $C>0$ such for all $t>0$, $p\geq1$, $x\in M$,
\[ \mathbb E_{t,x,x}\big[|\mathcal Hol^\omega_{0,t}(W)-1_G|^p\big]\leq(Cpt)^p. \]

Equivalently, there exists some $\eta>0$ such that
\[ \sup_{x\in M}\sup_{t>0}\mathbb E_{t,x,x}\left[\exp\left(\frac\eta t\big|\mathcal Hol^\omega_{0,t}(W)-1_G\big|\right)\right]<\infty. \]
\end{theorem}

\begin{remark} 
From the proof we expect that we can allow manifolds $M$ that are not compact, provided the supremum with respect to $x$ be taken over a compact subset of $M$, i.e. $C$ becomes a locally bounded function of $x$. In fact it cannot be made uniform, because a time change is the same as a space dilation, so it is possible to artificially inflate the constant around a certain point by shrinking the metric in a similar fashion in a neighbourhood.

The argument that needs checking is the comparison between Brownian bridges and paths. In the proof, one of the crucial ingredients is that the law of $W_{|[0,t/2]}$ under the Brownian bridge measure with horizon $t$ is not only absolutely continuous with respect to that under the unpinned Brownian measure, but also that its density is bounded (locally) uniformly in time and in starting point. This is not true for non-compact manifolds, however our computations suggest that for every stopping time $\tau\leq t/2$ such that $W_{|[0,\tau]}$ takes values in some relatively open set $U$, we should have
\[ \mathbb E_{t,x,x}[F(W_{\bullet\wedge\tau})]\leq C_U\mathbb E_x[F(W_{\bullet\wedge\tau})] \]
for any $F$ measurable bounded, uniformly in $t\leq1$ and $x\in U$. This would imply the local version of the theorem.
\end{remark}

\begin{remark}
The power of $t$ is optimal. In the next remark, we give an example where
\[   \liminf_{t\to0}\mathbb E^M_{t,x,x}\big[\big|\big(\mathcal Hol^\omega_{0,t}(W)-\Id\big)/t\big|^p\big]
\geq \mathbb E^{\mathbb R^2}_{1,0,0}\left[|\mathcal A(W)|^p\right]>0, \]
where $\mathcal A(W)$ is the Lévy area of a planar Brownian motion. This situation should be considered the generic one: we expect the distribution of $(\mathcal Hol^\omega_{0,t}(W)-1_G)/t$ under $\mathbb P_{t,x,x}$ to admit a non-trivial limit as $t\to0$, which only depends on the curvature of $\omega$ at $x$.
\end{remark}

\begin{remark}
In the typical case, for $\eta>0$ large enough,
\[ \lim_{t\to0}\mathbb E_{t,x,x}\big[\exp\big(\eta\big|\mathcal Hol^\omega_{0,t}(W)-1_G\big|/t\big)\big]=\infty, \]
so the power of $p^p$ is optimal in the first formulation; in other words, the tail is indeed exponential. Consider for instance the connection
\[ \nabla:s\mapsto\mathrm ds + i\frac\alpha2s(x\mathrm dy-y\mathrm dx) \]
on the trivial complex bundle over the torus $\mathbb{R}^2/\mathbb{Z}^2$, $\alpha>0$ fixed. We find
\[ \mathcal Hol^\nabla_{0,t}(W) = \exp\left(-i\alpha\mathcal A(W)\right) \]
whenever the diameter of $W$ is small enough (which by Lemma \ref{lem:exitball} becomes the overwhelmingly likely scenario), where $\mathcal A(W)$ is the Lévy area of the loop. By scale invariance, we get
\[   \mathbb E^{\mathbb R^2/\mathbb Z^2}_{t,x,x}\big[\exp\big(\eta\big|\mathcal Hol^\nabla_{0,t}(W)-\Id\big|/t\big)\big]
\geq \mathbb E^{\mathbb R^2}_{1,0,0}\left[\exp\left(\eta\left|e^{i\alpha t\mathcal A(W)}-1\right|/t\right)\mathbf1_{\operatorname{diam}(W)\leq 1/\sqrt t}\right],  \]
and by Fatou's lemma
\[   \liminf_{t\to0}\mathbb E^{\mathbb R^2/\mathbb Z^2}_{t,x,x}\big[\exp\big(\eta\big|\mathcal Hol^\nabla_{0,t}(W)-\Id\big|/t\big)\big]
\geq \mathbb E^{\mathbb R^2}_{1,0,0}\left[\exp\left(\eta\alpha|\mathcal A(W)|\right)\right]. \]
The Lévy area was proved by Lévy \cite{Levystochasticarea} to have the explicit Laplace transform
\[ \mathbb E^{\mathbb R^2}_{1,0,0}[\exp(\lambda\mathcal A(W))]=\frac\lambda{2\sin(\lambda/2)} \qquad\text{ for $|\lambda|<2\pi$.}
\]
Since the right-hand side has a pole at $\lambda=2\pi$, for $\eta\geq 2\pi/\alpha$ the lower bound above is infinite, and the corresponding exponential moments of the normalized holonomy are unbounded.

We stress that this is actually the generic situation, and that we only expect the holonomy to be better behaved if geometric conditions impose that the holonomies are trivial for \emph{every} curve (for instance when the connection is flat on a simply connected manifold like the spheres). Even if the connection is flat but not trivial, then for some free homotopy class $C$ of curves for which the holonomy $\mathcal H$ (independent of the curve) is not the identity,
\begin{align*}   \mathbb E_{t,x,x}\big[\exp\big(\eta\big|\mathcal Hol^\omega_{0,t}(W)-1_G\big|/t\big)\big]
&\geq \mathbb P_{t,x,x}\big(X\in C\big)\exp(\eta|\mathcal H-1_G|/t) \\
&\geq \exp((\eta|\mathcal H-1_G|-\delta)/t) \end{align*}
for $\delta>0$ small enough by the large deviation principle for Brownian bridges (a result due to Hsu \cite{HsuLDP}), which is clearly unbounded as $t\to0$, provided $\eta$ is chosen large enough.
\end{remark}

\begin{proof}[Proof of Theorem \ref{th:momentsoftheholonomy}]
We show the first estimate, and the second one follows by summation. Let $E$ be the vector bundle adjoint to $P$ via the representation of $G$ on $\mathbb R^n$, and let $\nabla$ be the connection induced on it by $\omega$. The holonomy of $\nabla$ over a closed curve is precisely the action of the holonomy of $\omega$ over the same curve, so
\[ \big|\mathcal Hol^\omega_{0,t}(W)-1_G\big|
 = \big|\mathcal Hol^\nabla_{0,t}(W)-\operatorname{Id}\big|. \]
We know that the laws of $W_{|[0,\frac{t}{2}]}$ under $\mathbb P_{t,x,x}$ and $\mathbb P_{x}$ are absolutely continuous with respect to each other, and in fact the density is uniformly bounded. Indeed, it is given explicitly \cite[Equation (5.4.2)]{Hsu} by  
\[ 
\frac{(W_{|[0,\frac{t}{2}]})_*\mathrm d\mathbb P_{t,x,x}}{(W_{|[0,\frac{t}{2}]})_*\mathrm d\mathbb P_x} =\frac{p_\frac{t}{2}(W_{\frac{t}{2}},x ) }{p_t(x,x)}.
\]
It follows from the semigroup property that for all $x\in M$, $p_t(x,x)=\| p_{\frac{t}{2}}(x, \cdot)\|_{L^2}^2$, and using the triangle inequality we deduce that for all $x,y\in M$, 
\[ 
p_t(x,y)= \int  p_{\frac{t}{2}}(x, z) p_{\frac{t}{2}}(z,y)\d z \leq \| p_{\frac{t}{2}}(x, \cdot)\|_{L^2}
\| p_{\frac{t}{2}}(y, \cdot)\|_{L^2}=\sqrt{p_{t}(x,x)p_{t}(y,y)}. 
\]
Using the heat kernel expansion \eqref{eq:Kexpansion} for $p_t$ at the first order, we deduce that there exist $t_0>0, C<\infty$ such that for all $t\leq t_0$ and all $x,y\in M$,  
\[ p_{\frac{t}{2}}(y,y) \leq C p_t(x,x). \]
All together, we have 
\[ 
\sup_{t\leq t_0, x\in M} \frac{(W_{|[0,t/2]})_*\mathrm d\mathbb P_{t,x,x}}{(W_{|[0,t/2]})_*\mathrm d\mathbb P_x}
\leq \sup_{t\leq t_0, x, y\in M} \frac{ p_{\frac{t}{2}}(x,y)}{p_t(x,x)}\leq\sup_{t\leq t_0, x, y\in M}  \frac{ \sqrt{p_{\frac{t}{2}}(x,x)p_{\frac{t}{2}}(y,y)}    }{p_t(x,x)}\leq C. 
\]

Thus, 
\[ |\rho_{t,x}|_\infty\coloneqq \sup\frac{(W_{|[0,t/2]})_*\mathrm d\mathbb P_{t,x,x}}{(W_{|[0,t/2]})_*\mathrm d\mathbb P_x}<\infty, \quad \mbox{and}\quad 
|\rho|_\infty\coloneqq \sup_{t\in[0,t_0],x\in M}|\rho_{t,x}|_\infty<\infty.
\]

Let us come back to Theorem \ref{th:momentsoftheholonomy}. Since $M$ is compact, it is enough to show that there exists an open cover by subspaces $U$ such that the estimate holds for all $x\in U$, with a constant $C$ depending on $U$. Let $\phi:U\to B_0(2)$ be a chart from an open set $U$ of $M$. The preimage $\phi^{-1}B_0(1/2)$ is open, and the collection of all these sets for $\phi$ ranging over all such charts covers $M$, so it will suffice to show the estimate in the open set.

Let $x\in M$ be such that $|\phi(x)|<1/2$, and let us find upper bounds that depend on $x$ only through $\phi$ and $U$. Let $\tau$ be the first time for which $W$ exits $\phi^{-1}B_0(1)$. We consider first the event $\{\tau\leq t\}$, and assume $t\leq t_0$ during the proof, which is enough to conclude since $|\mathcal Hol^\omega_{0,t}(W) -1_G|$ is globally bounded by $2$. We know by symmetry and the exit time estimate that
\begin{align*}
      \mathbb E_{t,x,x}\big[|\mathcal Hol^\nabla_{0,t}(W)-\mathrm{Id}|^p\mathbf1_{\tau\leq t}\big]
&\leq 2\cdot2^p\mathbb P_{t,x,x}(\tau\leq t/2) \\
&\leq 2|\rho|_\infty2^p\mathbb P_x(\tau\leq t/2) \\
&\leq C_0|\rho|_\infty2^{p+1}\exp(-2/(C_0t))
\end{align*}
for all $t>0$, $p\geq1$, and some $C_0>0$ independent of $x$. Since
\[ \exp(-2/(C_0t))\leq \Gamma(p+1)(C_0t/2)^p \]
for all $p,t>0$, we find
\[   \mathbb E_{t,x,x}\big[|\mathcal Hol^\nabla_{0,t}(W)-\mathrm{Id}|^p\mathbf1_{\tau\leq t}\big]
\leq \frac{2C_0|\rho|_\infty\Gamma(p+1)}{p^p}(C_0pt)^p, \]
which satisfies the expected bound since the first factor converges to zero exponentially fast as $p$ goes to infinity.

Let $A$ be a section of the restriction $(T^*M\otimes P_{\mathfrak g})_{|U}$ of the bundle to the open set $U$, satisfying the conditions of Lemma \ref{lem:flatmodification}, and write $\nabla'$ for $\nabla+(A_x\cdot{})$; recall that $\nabla'$ is flat by definition, and since $U$ is simply-connected its holonomy along any closed curve is the identity. On the event $\{t<\tau\}$, using the fact that the holonomy from $s$ to $t$ is adjoint to that from $t$ to $s$, hence has the same norm, we know that
\begin{align*}
      \big|\mathcal Hol^\nabla_{0,t}(W)-\Id\big|
&=    \big|\mathcal Hol^{\nabla'}_{0,t}(W)-\mathcal Hol^\nabla_{0,t}(W)\big|\\
&\leq \big|\mathcal Hol^{\nabla'}_{t/2,t}(W)\circ\big(\mathcal Hol^{\nabla'}_{0,t/2}(W)-\mathcal Hol^\nabla_{0,t/2}(W)\big)\big| \\
&     \qquad + \big|\big(\mathcal Hol^{\nabla'}_{t/2,t}(W)-\mathcal Hol^\nabla_{t/2,t}(W)\big)\circ\mathcal Hol^\nabla_{0,t/2}(W)\big| \\
&\leq \big|\mathcal Hol^{\nabla'}_{0,t/2}(W)-\mathcal Hol^\nabla_{0,t/2}(W)\big|  + \big|\mathcal Hol^{\nabla'}_{t/2,t}(W)-\mathcal Hol^\nabla_{t/2,t}(W)\big| \\
&=    \big|\mathcal Hol^{\nabla'}_{0,t/2}(W)-\mathcal Hol^\nabla_{0,t/2}(W)\big|  + \big|\mathcal Hol^{\nabla'}_{t,t/2}(W)-\mathcal Hol^\nabla_{t,t/2}(W)\big|.
\end{align*}
By symmetrisation and the density argument again,
\[ \mathbb E_{t,x,x}\big[|\mathcal Hol^{\nabla'}_{0,t}(W)-\mathcal Hol^\nabla_{0,t}(W)|^p\mathbf1_{t<\tau}\big]
\leq 2^p|\rho|_\infty\mathbb E_x\big[|\mathcal Hol^{\nabla'}_{0,t/2}(W)-\mathcal Hol^\nabla_{0,t/2}(W)|^p\mathbf1_{t/2<\tau}\big]. \]
It will be enough to show
\[   \mathbb E_x\big[|\mathcal Hol^{\nabla'}_{0,t\wedge\tau}(W)-\mathcal Hol^\nabla_{0,t\wedge\tau}(W)|^p\big]
\leq (Cpt)^p \]
for some $C>0$, and for all $t>0$, $p\geq1$. But the expectation depends only on the part of $M$ within the open set $\phi^{-1}B_0(1)$, so we can identify $U$ with $B_0(2)$, push all our objects to it, and recover the setting of the beginning of the section. This gives precisely the hypotheses of Lemma \ref{lem:topologylessholonomyestimate}, and we conclude by noting again that the control does not depend on the position of $x$ within $\phi^{-1}B_0(1/2)$.
\end{proof}

This concludes the proof of the first inequality in Proposition \ref{prop:tech}. The other ones follow from an easy derivation.

\begin{proof}[Proof of Proposition \ref{prop:tech}]
The first inequality is precisely Theorem \ref{th:momentsoftheholonomy}. The fact we can replace $\mathbb{E}_{t,x,x}$ in it is trivial, and the third inequality follows from the first by triangle inequality. We now prove the second inequality, 
\[ \sup_{x\in M}  \big| \mathbb{E}_{t,x,x}^m\big[\operatorname{Id}_{E_x}-\mathcal{H}ol^{\nabla}_{0,t}(W)) \big]\big| \leq Ct^2, \]

We claim that the main contribution to $\operatorname{Id}_{E_x}-\mathcal Hol^\nabla_{0,t}(W) $ is symmetric in distribution, so taking the absolute value in the expectation would lead to a worse (and indeed insufficient) bound. Instead, we get rid of this term by a symmetrisation argument, and bound the next term in the expansion. 
Since the Brownian bridge is invariant under time reversal, we have
\[ \mathbb E^m_{t,x,x}\big[\mathcal Hol^\nabla_{0,t}(W)\big] = \mathbb E_{t,x,x}^m\big[ \mathcal Hol^\nabla_{t,0}(W)\big] = \mathbb E_{t,x,x}^m\big[ \mathcal Hol^\nabla_{0,t}(W)^*\big]. \]
Noticing that for every isometry $U$,
\begin{equation} \big| \operatorname{Id}-\frac12(U+U^*)  \big|= |\frac12  U^*(\operatorname{Id}-U)^2 |\leq\frac12  |\operatorname{Id}-U|^2, \label{eq:sym}\end{equation}
we see that
\begin{align*}
      \big|\mathbb{E}^m_{t,x,x}\big[ 1-\mathcal{H}ol^{\nabla}_{0,t}(W) \big]\big|
&=       \big|\mathbb{E}^m_{t,x,x}\big[ \operatorname{Id}_{E_x}-\frac12(\mathcal{H}ol^{\nabla}_{0,t}(W)-\mathcal{H}ol^{\nabla}_{0,t}(W)^*)  \big]\big|\\
&\leq \mathbb{E}^m_{t,x,x}\big[ \big|\operatorname{Id}_{E_x}-\frac12(\mathcal{H}ol^{\nabla}_{0,t}(W)-\mathcal{H}ol^{\nabla}_{0,t}(W)^*) \big| \big]\\
&\leq \frac12\cdot\mathbb E_{t,x,x}\Big[\big|\operatorname{Id}_{E_x}-\mathcal Hol^\nabla_{0,t}(W)\big|^2\Big] \\
&\leq C(E,\nabla) t^2, \qquad \mbox{by Theorem \ref{th:momentsoftheholonomy}. }
\end{align*}
We easily deduce the fourth inequality: since the trace commutes with the expectation, we have 
\begin{align*}
\big|\mathbb{E}_{t,x,x}^m\big[ \chi(W)\big]\big|
&=\big|\mathbb{E}_{t,x,x}^m\big[ \tr\big(\operatorname{Id}_{E_{1,x}}- \mathcal Hol^{\nabla_1}_{0,t}(W) \big)     \big] - \mathbb{E}_{t,x,x}^m\big[ \tr\big(\operatorname{Id}_{E_{0,x}}-\mathcal  Hol^{\nabla_0}_{0,t}(W) \big)     \big]  \big|\\
&=\big|\tr\big(\mathbb{E}_{t,x,x}^m\big[ \operatorname{Id}_{E_{1,x}}- \mathcal Hol^{\nabla_1}_{0,t}(W)      \big]\big) - \tr\big(\mathbb{E}_{t,x,x}^m\big[ \operatorname{Id}_{E_{0,x}}- \mathcal Hol^{\nabla_0}_{0,t}(W)    \big] \big) \big|\\
&\leq \big|\mathbb{E}_{t,x,x}^m\big[ \operatorname{Id}_{E_{1,x}}- \mathcal Hol^{\nabla_1}_{0,t}(W)      \big]\big|
+\big|\mathbb{E}_{t,x,x}^m\big[ \operatorname{Id}_{E_{0,x}}- \mathcal Hol^{\nabla_0}_{0,t}(W)      \big]\big|\\
&\leq (C(E_1,\nabla_1)+C(E_0,\nabla_0) ) t^2,
\end{align*}
which concludes the proof.
\end{proof}

\section{Covariant Symanzik Identities}
\label{sec:Symanzik}
In this section, $(M,m)$ is assumed to be a manifold of dimension $d\in\{2,3\}$ without boundary and with non-vanishing non-negative mass. 
As opposed to the other sections, the connection $\nabla$ on $E$, rather than being fixed, is a random variable under some arbitrary probability distribution $\mathbf{P}$, with expectation $\mathbf{E}$, on the space of smooth connections. The interested readers may choose $E$ and its metric to be random as well by conditioning on the bundle, since there are only countably many of them up to isomorphism.%
\footnote{We could not find a reference for this precise result. Because isomorphic metric bundles are in fact isometric, it is enough to count isomorphism classes of bundles without metrics. Rank $n$ vector bundles over $M$ are classified by smooth maps $M\to G_{n,N}$ to the so-called Grassmannian manifold of $n$-planes in $\mathbb R^N$, $N=n+d$, up to some relation which makes the equivalence classes open in the continuous topology. One can then consider countable dense subsets of $\mathcal C^\infty(M,G_{n,N})$. For more details, see e.g. \cite{MilnorCC}.}
Because we will deal with several probability spaces, the expectation with respect to the (massive) Brownian loop soup of intensity $\alpha\geq0$ will be written $\mathbb{E}^\mathscr{L}_\alpha$ in this section, and $\mathbb{P}^\mathscr{L}_\alpha$ is the corresponding probability measure. The notations $\mathbb{E}$ will always be reserved for the expectation with respect to all the random objects into play.

\subsection{Regularity of holonomies as functions of the connection}
The space $\mathcal{E}_0$ of smooth connections on $E$ is an affine space modelled over the space of sections of $T^*M\otimes\End(E)$. It thus inherits any translation-invariant topology on that bundle, and we endow it with the $\mathcal{C}^{2+\varepsilon}$ topology, where $\varepsilon>0$ is arbitrary.%
\footnote{We are talking here about the Whitney topology. In our case, the weak and strong variants coincide, and it is the initial topology for the collection of restrictions $\nabla\mapsto\nabla_{|U}\in\mathcal C^{2+\varepsilon}(T^*V\otimes\End(V\times\mathbb R^n))$ for $U\simeq V\subset\mathbb R^d$ a trivializing open set in $M$.}
The corresponding $\sigma$-algebra coincides with the one induced by the uniform topology, and is also the same as the one induced by the smooth topology. Let $\mathcal{E}\subset \mathcal{E}_0$ be the closed affine subspace of \emph{metric} connections (modelled over the space of sections of $T^*M\otimes\mathfrak{o}(E)\subset T^*M\otimes\End(E)$). The product space $\mathcal{E}\times\mathscr{L}$ is endowed with the product $\sigma$-algebra.

\begin{lemma}
    \label{le:HolCont}
    Let $W$ be either a Brownian motion or a Brownian bridge. The random functional
    \[ 
    \nabla\in \mathcal{E}\mapsto  \mathcal Hol^\nabla(W)
    \]
    admits a modification which is almost surely continuous with respect to the  $\mathcal{C}^{2+\epsilon}$ topology.
\end{lemma}
\begin{proof}
    By cutting the path at the middle, the holonomy over a bridge decomposes as a product of two elements.
    Both of those are the holonomy over a path whose law is absolutely continuous with respect to a Brownian motion, so it suffices to prove the result for $W$ a Brownian motion.
    
    Recall that $W$ can be written as the solution of some stochastic differential equation. Since $\mathcal Hol^\nabla(W)$ solves a Stratonovich differential equation driven by $W$, namely \eqref{eq:SDEHol}, the composition gives a description of the holonomy as the solution to some Stratonovich differential equation driven by a standard $d$-dimensional Brownian motion $w$, namely \eqref{eq:SDEframe} where $U^E_s = \mathcal Hol^\nabla_{0,st}(W)$. In fact, we can invert this equation to define $w$ uniquely from $W$, up to an inessential choice of isometry $u:\mathbb R^d\to T_xM$.
    
    The result follows from the continuity of the solution map with respect to the coefficients, up to modification, which is true in a general framework and we prove it using rough paths techniques.  Considering the Stratonovich path lift $\mathbf w$ of $w$, it holds that for all $\nabla\in \mathcal{E}$, almost surely, $\mathcal Hol^\nabla(W)$ is equal to $\mathcal{H}ol^\nabla(\mathbf w)$, the solution to the same equation but interpreted as a rough differential equation driven by $\mathbf w$, rather than as a Stratonovich SDE. Almost surely, the RDE solution $\mathcal{H}ol^\nabla(\mathbf w)$ is defined for all $\nabla$ simultaneously. Furthermore, the coefficients depend continuously on $\nabla$, when we consider the $\mathcal C^{2+\varepsilon}$ topology on both sides. Theorem 3 in \cite{CoutinFrizVictoir} precisely states that the solution of an RDE with a driver in $\mathcal{C}^{1/p}$ (in the rough path sense) depends continuously on the coefficients, when these are endowed with the topology of $\mathcal{C}^{p+\varepsilon'}$, for any $\varepsilon'>0$. Since the Brownian rough path has such a regularity for all $p>2$, it suffices to apply the result to $p=2+\varepsilon/2$ to conclude.
\end{proof}
From now on and with the previous lemma in mind, we write $\mathcal Hol^\nabla(W)$ for this continuous modification.
\begin{corollary}
    The map 
        \[ 
    (\nabla, W) \mapsto  \mathcal Hol^\nabla(W)
    \]
    from $\mathcal E\times \mathscr{L}$ is measurable. 

\end{corollary}
\begin{proof}
    This map is measurable in $W$ for all $\nabla$ and continuous in $\nabla$ for all $W$, which makes it a Carathéodory function. Since the $\mathcal{C}^{2+\epsilon}$ topology on the space $\mathcal{C}^\infty$ is metrisable separable, the conclusion follows from generic arguments, see e.g. \cite[Theorem 4.51]{Hitchhiker}.
\end{proof}
\begin{corollary}
    The map 
    \[
    \nabla\mapsto \mathbb{E}^{\mathscr L}_\alpha\bigg[ \prod_{\ell \in \mathcal{L}} \tr(\mathcal Hol^\nabla(\ell))\bigg]
    \]
    is measurable. Recall that the product is not absolutely convergent.
\end{corollary}
\begin{proof}
    By lemma \ref{le:HolCont}, the map $\nabla\mapsto \mathcal Hol^\nabla (W)$ is continuous. By dominated convergence (recall that $|\tr(Hol^\nabla (W))|\leq 1$ and $p_t(x,x)$ decays exponentially), we deduce that for any $\delta>0$, the map 
    \[ 
    f:\nabla \mapsto \int_\delta^\infty \int_M \frac{p_t(x,x)}{t} \mathbb{E}_{t,x,x}^m\big[  \tr(\mathcal Hol^\nabla (W))\big] \d x \d t
    \]
    is also continuous. 
    
    For any $\delta>0$, and setting $\lambda=\int_\delta^\infty \int_M \frac{p_t(x,x)}{t}\mathbb{E}^m_{t,x,x}[ \mathbf{1} ] \d x \d t<\infty$, we have
     \begin{align*}
    \mathbb{E}^{\mathscr L}_\alpha\bigg[ \prod_{\ell \in \mathcal{L}_\delta} \tr(\mathcal Hol^\nabla(\ell))\bigg]
    &=e^{-\alpha\lambda} \sum_{k=0}^\infty \frac{(\alpha\lambda)^k }{k!} \mathbb{E}^{\mathscr L}_\alpha\bigg[ \prod_{\ell \in \mathcal{L}_\delta} \tr(\mathcal Hol^\nabla(\ell) )\bigg|  \#\mathcal{L}_\delta=k \bigg]\\
    &=e^{-\alpha\lambda}
    \sum_{k=0}^\infty \frac{(\alpha f(\nabla))^k}{k!}=e^{\alpha(f(\nabla)-\lambda)},
     \end{align*}
    which is continuous in $\nabla$ since $f$ is. We conclude by convergence of the product in the $L^1$ sense as $\delta\to0$, and the fact that pointwise limits of continuous functions are measurable. 
\end{proof}
\begin{remark}
    The map $\nabla\mapsto\mathbb{E}^{\mathscr L}_\alpha[\prod_{\ell\in\mathcal L}\operatorname{tr}(\mathcal Hol^\nabla(\ell))]$ is in fact not only measurable but continuous. The proof we just presented readily extends, the missing argument is to show that
    \[ \sup_{\nabla\in B}\int_0^\delta \int_M \frac{p_t(x,x)}{t} \big|\mathbb{E}_{t,x,x}^m\big[  \tr(\mathcal Hol^\nabla (W))-1\big]\big| \d x \d t \]
    can be made arbitrarily small by taking $\delta$ sufficiently small, where $B$ is an arbitrary bounded set in $\mathcal C^{2+\varepsilon}$. This, in turn, follows from the fact that the technical estimation Proposition \ref{prop:tech} holds locally uniformly over $\nabla$. This can be tracked in the proof of Proposition \ref{prop:tech}, where in fact we could even use the topology $\mathcal{C}^1$.
\end{remark}

\subsection{Gaussian free fields twisted by a connection}

Informally, a (massive) Gaussian free field with values in $E$ twisted by $\nabla$, which implicitly depends on $m$, $g$ and the metric on $E$, is a centred Gaussian random section of $E$ with respect to the quadratic energy
\[
\frac{\| \nabla \Phi \|^2 }{2}+\langle \Phi, m \Phi\rangle.
\]
It is often written in the physics fashion, namely
\[
\frac{1}{Z_\nabla} \exp\Big( -\frac{\| \nabla \Phi \|^2 }{4}-\frac{\langle \Phi, m \Phi\rangle }{2} \Big) \mathcal{D}\Phi
\]
for $\mathcal D\Phi$ the putative volume element in the space of all sections.
Mathematically, $\Phi$ is a random \emph{distribution} in the sense of Schwartz
such that $\Phi(s)$ is Gaussian centred for all $s\in\mathcal C^\infty(M,E)$, and such that
    \[ 
    \mathbb{E}^\mathrm{GFF}_\nabla[ |\Phi(s)|^2]= \int_{M^2} \langle s(x) , G_\nabla(x,y) s(y) \rangle \d x \d y.
    \]
Recall that $G_\nabla$ is the Green kernel for $L_\nabla$. For $\nabla\in\mathcal E$ fixed, the probability measure $\mathbb{P}^\mathrm{GFF}_\nabla$ is known to exist on the space $\mathcal D'(M,E)$ of Schwartz distributions with values in $E$. The underlying $\sigma$-algebra is the Borel algebra for the strong dual topology, and is generated by the cylindrical events $\{\Phi(s)\in A\}$. In fact, $\Phi$ almost surely lies in the Sobolev spaces $H^{1-d/2-\varepsilon}(M,E)$ for all $\varepsilon>0$, and we could carry our analysis using this measurable space instead. See for example \cite{CoursBerestycki,CoursWernerPowell} for surveys on the scalar Gaussian free field, and \cite{Schwartz} for general considerations on $\sigma$-algebras and measures on spaces of distributions.

\begin{lemma}
    \label{le:twistedcovariance}
    For a smooth $\nabla$ and $\Phi$ a Gaussian free field twisted by $\nabla$, for $s_1,s_2$ two smooth sections of $E$,
    the covariance rewrites as
    \[\mathbb{E}^\mathrm{GFF}_\nabla\big[\Phi(s_1)\overline{\Phi(s_2)}\big]=\int_0^\infty \int_{M^2} p_t(x,y) \mathbb{E}^m_{t,x,y}[\langle s_1(x),\mathcal Hol^\nabla_{t,0}(W)s_2(y)\rangle ]\mathrm dx\,\mathrm dy\,\mathrm dt. \]
    In particular, it is continuous in $\nabla$, everything else being fixed.

    Moreover, for any smooth sections $s_1,\dots , s_k\in \mathcal{C}^\infty(M,E)$ and continuous bounded functions $f_1,\dots f_k: \mathbb C\to \mathbb{R}$, the functional
    \[ \nabla\mapsto\mathbb{E}^\mathrm{GFF}_\nabla\big[f_1(\Phi(s_1)) \dots f_k(\Phi(s_k))  \big] \]
    is also continuous.
\end{lemma}
\begin{proof}
    The first formula follows from unwrapping the definition of $\Phi$, then the definition of $G_\nabla$, and then applying Theorem \ref{th:K=pE}.

    By lemma \ref{le:HolCont}, for all $(t,x,y)\in(0,\infty)\times \mathbb{R}^2$, $\mathbb{P}_{t,x,y}^m$-almost surely, $\nabla\mapsto \mathcal Hol^\nabla(W)$ is continuous in $\nabla$. By the dominated convergence theorem, it follows that 
    the right-hand side in the lemma is also continuous in $\nabla$ (still in the $\mathcal{C}^{2+\epsilon}$ topology). 
    The domination comes from the fact that
    \begin{align*} \int_0^1 \int_{M^2} p_t(x,y) \mathrm dx\,\mathrm dy\,\mathrm dt&\leq \operatorname{vol}(M), &
    \int_1^\infty \int_{M^2} p_t(x,y) \mathbb{E}^m_{t,x,y}[\mathbf1 ]\mathrm dx\,\mathrm dy\,\mathrm dt&<\infty.\end{align*}
    The second bound comes from the spectral gap for the semigroup $(e^{-tL_{\mathrm d}})_{t\geq0}$ associated to the massive Laplace-Beltrami operator $L_{\mathrm d}=-\frac{1}{2}\Delta+m$, for instance writing
   \[  \int_1^\infty \int_{M^2} p_t(x,y) \mathbb{E}^m_{t,x,y}[\mathbf1 ]\mathrm dx\,\mathrm dy\,\mathrm dt
    = \int_1^\infty\big\langle\mathbf1,e^{-tL_{\mathrm d}}\mathbf1\big\rangle\mathrm dt
    \leq \operatorname{vol}(M)\int_1^\infty e^{-t \mu_0}\mathrm dt, \]
    where $\mu_0>0$ is the smallest eigenvalue of $L_{\mathrm d}$.
    
    As for the second, since the topology on $\mathcal E$ is metrisable, we need only prove sequential continuity. We consider a converging sequence $(\nabla_n)_{n\geq0}$ of connections, all other objects being fixed. Under $\mathbb P^\mathrm{GFF}_{\nabla_n}$, the vector $(\Phi(s_1),\ldots,\Phi(s_k))$ is centred Gaussian, with covariance given by the first formula. Since this covariance converges to the covariance of the limit by the continuity we have just shown,   the vectors converge in distribution to the limiting Gaussian vector, and by composition the variables $f_1(\Phi(s_1))\cdots f_k(\Phi(s_k))$, seen as variables under the sequence of probability measures $\mathbb P^\mathrm{GFF}_{\nabla_n}$, converge in distribution to the same variable under the limit measure. Since this variable is bounded, it implies convergence in moment, which is the expected conclusion.
\end{proof}

\subsection{Mixture of Gaussian free fields depending on a random connection}
\label{ssec:minimalcouplingGFF}

We want to define a new probability measure $\mathbb{P}^\mathrm{GFF}_\mathbf{P}$ on the set of couples $(\nabla,\Phi)$, \emph{formally} defined by 
\[
\d \mathbb{P}^\mathrm{GFF}_\mathbf{P}(\nabla,\Phi) = \frac{1}{Z} \exp\Big( -\frac{\| \nabla \Phi \|^2 }{4}-\frac{\langle \Phi, m\Phi\rangle\rangle }{2} \Big) \mathcal{D}\Phi \d \mathbf{P}(\nabla),
\]
where $\mathcal{D}\Phi$ is formally a Lebesgue measure on the vector space of sections of $E$, and $Z$ is a normalisation constant which crucially \emph{should not depend upon $\nabla$}.

With formal computations, we can rewrite this as follows:
\[
\d \mathbb{P}^\mathrm{GFF}_\mathbf{P}(\nabla,\Phi) = \frac{\Det(L_\nabla/2\pi)^{-\frac{1}{2}} }{Z} \frac{1}{\Det( L_\nabla/2\pi)^{-\frac{1}{2}} }\exp\Big( -\frac{\langle \Phi, L_\nabla\Phi\rangle^2 }{2} \Big) \mathcal{D}\Phi \d \mathbf{P}(\nabla).
\]
For a given $\nabla$, the term 
\[\frac{1}{\Det(L_\nabla/2\pi)^{-\frac{1}{2}} }\exp\Big( -\frac{\langle \Phi, L_\nabla\Phi\rangle^2 }{2} \Big) \mathcal{D}\Phi \]
should then be interpreted as the massive Gaussian free field twisted by $\nabla$. It is thus already normalised to have, $\nabla$-almost surely, a unit total mass, and we deduce that the weight
\[ \frac{\Det(L_\nabla/2\pi)^{-\frac{1}{2}} }Z \]
should be equal to
\[ \frac{\Det(L_\nabla)^{-\frac{1}{2}} }{\mathbf E'[\Det(L_{\nabla'})^{-\frac{1}{2}}]},\]
where $\nabla'$ is a random connection under $\mathbf P'$ whose distribution is $\mathbf P$, which we think of as an independent copy of $\nabla$.
As the expectation of an exponential, it is clear that the expectation ${\mathbf E'[\Det(L_{\nabla'})^{-\frac{1}{2}}]}$ is well defined and positive, although it might be infinite. This last case is ruled out by our representation as product of loops, see equation \eqref{eq:diamagnetic}.
Although this computation started with a formal object which has no rigorous definition, we  end up with a decomposition of it as a product of perfectly well-defined quantities which we are able to use as a definition.

\begin{prop}
    \label{prop:measurability}
    The annealed probability measure on $\mathcal{E}\times \mathcal{C}^\infty(M,E)^*$ given by the formula 
    \[\int_{\mathcal{E}} \delta_\nabla\otimes \mathbb{P}^\mathrm{GFF}_\nabla   \d \mathbf{P}(\nabla)\]
    is well defined. We denote its volume element by $\d \mathbb{P}^\mathrm{GFF}_\nabla (\Phi)\d \mathbf{P}(\nabla)$.
\end{prop}

\begin{proof}
    It suffices to show that for any measurable bounded $f:\mathcal E\times\mathcal{C}^\infty(M,E)^*\to \mathbb{R}$, the functional $\int f(\nabla,\Phi) \d \mathbb{P}^\mathrm{GFF}_\nabla (\Phi)$ is measurable in $\nabla$. By a monotone class argument, we can reduce to functions of the form \[ (\nabla,\Phi)\mapsto\mathbf1_{\nabla\in A}\prod_{i=1}^kf_i(\Phi(s_i)), \]
    for $f_i:\mathbb R\to\mathbb R$ smooth with compact support and $A\subseteq\mathcal E$ measurable. This is an immediate consequence of the continuity of
    \[ \nabla\mapsto\mathbb E^\mathrm{GFF}_\nabla\bigg[\prod_{i=1}^kf_i(\Phi(s_i))\bigg], \]
    which was established in Lemma \ref{le:twistedcovariance}.
\end{proof}

\begin{definition}
\label{def:Symanzik}
    The probability measure $\mathbb{P}^\mathrm{GFF}_{\mathbf{P}}$ heuristically described at the beginning of the section is defined as
    \[
    \d  \mathbb{P}^\mathrm{GFF}_\mathbf{P}(\nabla, \Phi)=\frac{\exp\big(\frac{\zeta'_{\nabla}(0)}{2}\big)  } { \mathbf{E}'\left[\exp\big( \frac{\zeta'_{\nabla'}(0)}{2}  \big)\right]  }  \d \mathbb{P}^\mathrm{GFF}_\nabla (\Phi) \d \mathbf{P}(\nabla);\]
    equivalently,
    \[\d  \mathbb{P}^\mathrm{GFF}_\mathbf{P}(\nabla, \Phi)=\frac{  \mathbb{E}^\mathscr{L}_{n/2}[ \prod_{\ell\in \mathcal{L}} \tr  \mathcal Hol^\nabla(\ell)  ]    }{ \mathbf{E}'\otimes\mathbb{E}^\mathscr{L}_{n/2}[ \prod_{\ell\in \mathcal{L}}  \tr \mathcal Hol^{\nabla'}(\ell)  ] } \d \mathbb{P}^\mathrm{GFF}_\nabla (\Phi) \d \mathbf{P}(\nabla). \]
\end{definition}

The following theorem is known in the physics literature as a Symanzik formula. A discrete version is given by Lévy and Kassel in \cite{KasselLevy}, while the flat case is treated in \cite{LeJan}.

\begin{theorem}
    \label{th:symanzik}
Assume $\mathbb{K}=\mathbb{R}$.    Let $(s_1,\dots, s_{2k} )$ be smooth sections of $E$. Let $f$ be a measurable functional of $\nabla$. Then we have
    \begin{multline*} 
    \mathbb{E}^\mathrm{GFF}_\mathbf{P}[ f(\nabla)\Phi(s_1)\dots \Phi(s_{2k }) ]\\
    =
    \sum_{\pi} 
    \frac{ \mathbf{E}\Big[f(\nabla )  \mathbb{E}^{\mathscr{L}}_{n/2}[\mathcal{Z}_{\nabla,\mathcal{L}}]  \prod_{\{i,j\}\in \pi } \int_0^\infty \int_{M^2} \mathbb{E}_{t,x,y}^m [ \langle s_i(x), \mathcal Hol_{t,0}^\nabla(W) s_j(y) \rangle]\d x\d y \d t   \Big] }{ 
    \mathbf{E}\otimes\mathbb E^\mathscr{L}_{n/2}[\mathcal{Z}_{\nabla,\mathcal{L}}]
    },
    \end{multline*}
    where
    \[ 
    \mathcal{Z}_{\nabla,\mathcal{L}}= \prod_{\ell \in \mathcal{L}}  \tr(\mathcal Hol^\nabla(\ell))
    \]
    and where $\pi$ ranges over pairings of $\{1,\dots, 2k\}$ (i.e. there exist $i_1, \dots, i_k, j_1, \dots, j_k\in \{1,\dots, 2k\}$ such that $\{i_1,\dots,i_k,j_1,\dots,j_k\}=\{1,\dots, 2k\}$ and $\pi=\{ \{i_l,j_l\},1\leq l\leq k\}$ ).

    The case $\mathbb{K}=\mathbb{C}$ is identical, except the sum is over permutations of $\{1,\ldots,k\}$ rather than pairings of $\{1,\ldots,2k\}$:
    \begin{multline*} 
    \mathbb{E}^\mathrm{GFF}_\mathbf{P}\big[ f(\nabla)\Phi(s_1)\dots\Phi(s_k)\overline{\Phi(s_{k+1})\cdots \Phi(s_{2k })} \big]\\
    =
    \sum_{\sigma}
    \frac{ \mathbf{E}\Big[f(\nabla )  \mathbb{E}^{\mathscr{L}}_{n/2}[\mathcal{Z}_{\nabla,\mathcal{L}}]  \prod_{i=1}^k \int_0^\infty \int_{M^2} \mathbb{E}_{t,x,y}^m [ \langle s_i(x), \mathcal Hol_{t,0}^\nabla(W) s_{k+\sigma(i)}(y)\rangle ]\d x\d y \d t   \Big] }{ 
    \mathbf{E}\otimes\mathbb E^\mathscr{L}_{n/2}[\mathcal{Z}_{\nabla,\mathcal{L}}]
    }.
    \end{multline*}
\end{theorem}
\begin{proof}
    This now follows directly from the second representation formula for $\mathbb{P}^{\mathrm{GFF}}_{\mathbf{P}}$, Lemma \ref{le:twistedcovariance} for the second moment of a Gaussian free field twisted by $\nabla$, and Isserlis' theorem (also known as Wick's theorem) for higher moments of Gaussian random variables, see e.g. \cite[Theorem 1.36]{Janson}.
\end{proof}
Similar formulas, either as formal expression or in other settings, can be found for example in \cite[(79) and (80)]{KasselLevy} in a discrete setting, or in \cite[p. 134]{AlbeverioKusuoka} for the Abelian case over $\mathbb{R}^2$: from this second expression, by removing the local time contributions that comes from the $\Phi^4$ term,  we obtain our expression in the case $\mathbf{P}$ is the Yang--Mills measure, except the partition function is expressed in term of the intensity of the loop soup measure rather than in term of the loop soup itself. See also \cite{AHKHK} for related expressions which also connect similar observables of random fields coupled to random gauge fields with holonomies. There, the curvature of the connection is somehow discretised, which makes the small time divergences a bit different---see e.g. the paragraph between (4.1) and (4.2) in \cite{AHKHK}.

\section{Conformal invariance in dimension 2}
\label{sec:inv}
We work under the usual assumption that $\partial M\neq\emptyset$ or $m\neq0$. In this section, we consider specifically the case $d=2$, and we compare the partition functions of different metrics. We write $Z_{g,\nabla,m}$ for $\exp(\zeta'_L(0)/2)$, where $L$ is defined by using the Riemannian metric $g$ on $M$ and the mass term $m$. Similarly, we denote by $\Lambda_{g,m}$ the Brownian loop soup measure defined using the metric $g$ and the mass $m$, and $\mathcal{L}_{g,m}$ for the corresponding Brownian loop soup, with intensity $\alpha$. We define the random variable
\[ \mathcal{Z}_{g,\nabla,m}\coloneqq \prod_{\ell \in \mathcal{L}_{g,m}}  \tr(\mathcal Hol_{0,t}^\nabla(\ell)), \]
which, by Proposition \ref{prop:loops} applied with $\nabla_1=\nabla$ and $\nabla_0=\d $ on the trivial line bundle over $M$, is well defined as an almost sure limit of finite products. 

Our goal is precisely to show that the ratio of partition functions does not depend on $g$ in a given conformal class.
\begin{theorem}
    \label{th:conformal}
    Let $d=2$, let $(E,h_E,\nabla)$ be a metric bundle of rank $n$ and with a compatible connection. Let also $g$ be a Riemannian metric on $M$, $f,m\in \mathcal{C}^\infty(M,\mathbb{R})$ with $m\geq0$, and recall that either $M$ has a non-empty boundary, or that $m$ is not identically vanishing.
    Then, in distribution,
    \[ 
    \mathcal{Z}_{g,\nabla, m} \overset{(d)}{=} \mathcal{Z}_{e^{2f}g,\nabla, e^{-2f} m}.
    \]
    In particular, if $(E',h_{E'},\nabla')$ is another such triple on $M$ of rank $n'$,
    \begin{equation}
    \label{eq:conformal}
    \frac{Z_{g,\nabla,m}^{1/n}   }{Z_{g,\nabla',m}^{1/n'} }=\frac{  \mathbb{E}^\mathscr{L}_{1/2}[ \mathcal{Z}_{g,\nabla, m}    ]   }{ \mathbb{E}^\mathscr{L}_{1/2}[\mathcal{Z}_{g,\nabla', m}    ]  }=\frac{  \mathbb{E}^\mathscr{L}_{1/2}[ \mathcal{Z}_{e^{2f}g,\nabla,e^{-2f} m}    ]   }{ \mathbb{E}^\mathscr{L}_{1/2}[\mathcal{Z}_{e^{2f}g,\nabla', e^{-2f}m}    ]  }
    =\frac{Z_{e^{2f}g,\nabla,e^{-2f}m}^{1/n}}{Z_{e^{2f}g,\nabla',e^{-2f}m}^{1/n'}}.
    \end{equation}
\end{theorem}
\begin{remark}
    There are two ways to use this formula: either to say that the ratio for different connections does not depend on the metric within a conformal class, or to say that the ratio for different metrics do not depend on the specific bundle or connection we consider.

    In particular, this can be taken as a \emph{definition} for the ratio $\frac{Z_{g,\nabla', m}^{1/n'}}{Z_{g,\nabla,m}^{1/n}  }$ even when $g$ is not smooth (but provided $\nabla$ and $\nabla'$ are) --- we specifically have in mind the case when $g$ is a Liouville quantum gravity metric. 
\end{remark}

Let $\pi(\ell)$ be the class of a loop $\ell$ modulo rerooting and orientation-preserving reparametrisation, 
i.e. $\pi(\ell)=\pi(\ell')$ if and only if there exists $\theta\in [0,t_{\ell'}]$ and an increasing homeomorphism $\phi: [0,t_\ell]\to [0,t_{\ell'}] $ such that for all $t\in [0,t_\ell]$, $\ell(t)=\ell'( \phi(t)+\theta \mod t_{\ell'}) $. For a set of loops $\mathcal{L}$, $\pi(\mathcal{L})\coloneqq \{\pi(\ell), \ell\in \mathcal{L}\}$. 
\begin{lemma}
\label{le:reparam}
    Let $f\in \mathcal{C}^\infty(M,\mathbb{R})$. Then, the measures $\pi_*(\Lambda_{g,m})$ and 
    $\pi_*(\Lambda_{e^{2f}g,e^{2f} m})$ are equal. Equivalently, the random sets $\pi(\mathcal{L}_{g,m})$ and $\pi(\mathcal{L}_{e^{2f}g,e^{-2f} m})$ are equal in distribution. 
\end{lemma}

This is well known and follows essentially from \cite{Werner}, we only sketch the proof.
Explicitly, we reparameterise any loop $\ell\in\mathcal{L}_{g,m}$ of $g$-quadratic variation $T_\ell$ into a loop $\hat{\ell}$ of $e^{2f}g$-quadratic variation duration $\hat T_\ell$. To do this, we first reroot the loop uniformly with respect to the $e^{2f}g$-quadratic variation, then use the parameterisation by $e^{2f}g$-quadratic variation. The collection $\hat{\mathcal L}$ of the loops so obtained is distributed as $\mathcal{L}_{e^{2f}g,e^{-2f} m}$.

\begin{proof}[Proof of Lemma \ref{le:reparam} (sketch)]
    First, we assume $m=0$. Let $\mathcal{L}_0=\mathcal{L}_{g,0}$.
    For  $\tilde{g}\in \{g,e^{2f}g\}$ and $\ell\in \mathcal{L}$, let $\langle \ell\rangle_{\tilde{g},t}$ be the $\tilde{g}$-quadratic variation of $\ell_{|[0,t]}$. 
    Let $T_\ell$ be the total $g$-quadratic variation of $\ell$, and $\hat{T}_\ell$ its total $e^{2f}g$-quadratic variation. 
    Almost surely, for all $\ell\in \mathcal{L}_0$ and $t\in [0,T_\ell]$, 
    \begin{equation}
    \label{eq:sandwich1}
    \langle \ell\rangle_{e^{2f}g,t }=\int_0^t e^{2f(\ell_s)} \d \langle \ell\rangle_{g,s}=\int_0^t e^{2f(\ell_s)} \d s.
    \end{equation}
    In particular, $t\in [0,T_\ell]\mapsto \langle \ell\rangle_{e^{2f}g,t }$ is monotonically and continuously increasing. Let $\tau$ be its inverse (defined from $[ 0,\hat{T}_\ell] $ to $[0,T_\ell]$), so that $\ell\circ \tau$ is parameterized by unit $e^{2f} g$-quadratic variation. We now choose a new root for $\ell\circ \tau$: let $\theta_\ell$ be distributed uniformly over $[0, \hat{T}_\ell ]$, independently from $\mathcal{L}_0$ conditionally on $\hat{T}_\ell$, and independently for each $\ell$,%
    \footnote{Rigorously, we could choose a Poisson process of $(\ell,u_\ell)$ over $\mathscr L\times[0,1]$ with the product measure $\alpha\Lambda_{g,0}\otimes \operatorname{leb} $, and set $\theta_\ell\coloneqq u_\ell\hat T_\ell$.} and define 
    \[
    \hat{\ell}:t\in [0,\hat{T}_\ell] \mapsto \ell\circ \tau(  t+ \theta_\ell \mod \hat{T}_\ell ), 
    \]
    where $t+ \theta_\ell \mod \hat{T}_\ell$ is the representative in $[0, \hat{T}_\ell)$.
    We set $\hat{\mathcal{L}}_0=\{ \hat{\ell}: \ell \in \mathcal{L}_0 \}$. 
    It is easily seen that $\pi(\ell)=\pi(\hat{\ell})$ for all $\ell \in \mathcal{L}_0$, and $\hat{\mathcal{L}}_0$ is in fact distributed exactly as $\mathcal{L}_{e^{2f}g, 0}$ (not only up to rerooting and reparameterisation); this follows from the computation in \cite{Werner}. 
    
    To deal with the massive case $m\neq0$, we now proceed as follows. For each $\ell \in \mathcal{L}_0$, let $B_\ell$ be a Bernoulli random variable equal to $1$ with probability $\exp(- \int_0^{T_\ell} m(\ell_t)\d t)$; once again, this is defined with as much independence as possible, so that $\{ (\ell,\theta_\ell, B_\ell)\}$ is a Poisson point process.
    Then, the subset $\{ \ell\in \mathcal{L}_0: B_\ell=1\} $ is distributed as $\mathcal{L}_{g,m}$, so we can effectively define $\mathcal{L}_{g,m}$ as this set. Thus, it suffices to show that $\hat{\mathcal{L}}_m\coloneqq \{ \hat{\ell}: \ell \in  \mathcal{L}_{g,m} \}$ is distributed as $\mathcal{L}_{e^{2f}g,e^{2f} m}$. By the same argument, if $(\hat{B}_{\hat{\ell}})_{\hat{\ell} \in \hat{\mathcal{L}}_0}$ are Bernoulli random variables such that for each $\hat{\ell} \in \hat{\mathcal{L}}_0$,  $\hat{B}_{\hat{\ell}}$ is equal to $1$ with probability 
    $\exp(- \int_0^{\hat{T}_{\ell} } \hat{m}(\hat{\ell}_t)\d t) $ (for some function $\hat{m})$, and so that $\{(\hat{\ell}, \hat{B}_{\hat{\ell}})\}$ is a Poisson point process, then 
    $\{ \hat{\ell}\in \hat{\mathcal{L}}_0: \hat{B}_{\hat{\ell}}=1 \} $ is distributed as $\mathcal{L}_{e^{2f}g,\hat{m}} $.
    Thus, it suffices to show that for all $\ell\in \mathcal{L}_0$,
    \[ 
    \int_0^{\hat{T}_{\ell}  }  (e^{-2f} m)(\hat{\ell}_t)\d t=\int_0^{T_\ell} m(\ell_t)\d t . 
    \]
    For any $\ell \in \mathcal{L}_0$, with the change of variable $u=\tau(t)$ (thus $\d u= e^{-2 f(\ell_u)} \d t$), we get
    \[
    \int_0^{\hat{T}_{\ell}  } (e^{-2f}m)(\hat{\ell}_t)\d t 
    = \int_0^{T_\ell } (e^{-2f} m )( \ell_u)    e^{2 f(\ell_u)} \d u 
    = \int_0^{T_\ell } m( \ell_u)    \d u,  
    \]
    which concludes the proof.
\end{proof}
We now set $\mathcal{L}=\mathcal{L}_{g,m}$ and $\hat{\mathcal L}=\hat{\mathcal L}_m$ from the proof above, and keep the notations $\hat{\ell}$, $T_\ell$ and $\hat{T}_\ell$. Notice that $\mathcal Hol^{\nabla}(\ell)$ is invariant by reparametrisation of $\ell$. Since rerooting $\ell$ has the effect of conjugating $\mathcal Hol^{\nabla}(\ell)$ and the trace is conjugation-invariant, 
 $\tr( \mathcal Hol^{\nabla}(\ell))$ is invariant by reparametrisation and rerooting of the loop $\ell$. In particular, $\tr(\mathcal Hol^{\nabla}(\ell) )=\tr(\mathcal Hol^{\nabla}(\hat{\ell}) ) $ for all $\ell\in \mathcal{L}$. 

Here one must be aware of a notational subtlety: although we have shown that $\pi(\mathcal{L})=\pi(\hat{\mathcal{L}}) $ and that  $\tr( \mathcal Hol^{\nabla}(\ell))$ do not depend on the representative $\ell$ in $\pi(\ell)$, we cannot deduce yet that 
\[ \prod_{\ell\in \mathcal{L} } \tr( \mathcal Hol^{\nabla}(\ell)) = \prod_{\ell\in \hat{\mathcal{L}} } \tr( \mathcal Hol^{\nabla}(\ell)). \]
This would be automatic if the products were absolutely convergent, or if they were defined as a limit of finite products in a way that would not depend on $g$ nor $m$, but this is not the case and the proof requires some more work. 

\begin{proof}[Proof of Theorem \ref{th:conformal}]
We let $\mathcal{L}_\delta$, as before, be the subset of $\mathcal{L}$ of loops with $g$-quadratic variation larger than $\delta$, $\hat{\mathcal{L}}_\delta$ the subset of $\hat{\mathcal{L}}$ of loops with $e^{2f}g$-quadratic variation larger than $\delta$. Moreover, we write
$Q_\delta=\prod_{\ell \in \mathcal{L}_\delta } \tr(\mathcal Hol^{\nabla}({\ell}) ) $ and
$\hat{Q}_\delta=\prod_{\hat{\ell} \in \hat{\mathcal{L}}_\delta } \tr(\mathcal Hol^{\nabla}(\hat{\ell}) ) $ for the corresponding products.

For all $\ell \in \mathcal{L}$, it follows directly from \eqref{eq:sandwich1} that
\[ 
  \hat{T}_\ell\geq T_\ell/ \| e^{-2f}\|_\infty ,
\]
where
$\| e^{-2f}\|_\infty$ is positive and finite by compactness. Thus, for all $\delta>0$, 
\[ 
\pi(\mathcal{L}_{\delta} ) \subseteq \pi(\hat{\mathcal{L}}_{\delta/\| e^{-2f}\|_\infty}).
\]
Since by definition $\mathcal{Z}_{g,\nabla,m}=\lim Q_\delta $ and $\mathcal{Z}_{e^{2f}g,\nabla,e^{-2f} m}=\lim \hat{Q}_\delta $, it suffices to show  the convergence of $\hat Q_{\delta/\| e^{-2f}\|_\infty}-Q_{\delta}$ toward $0$ in distribution.
This difference rewrites as
\begin{align*}
\hat Q_{\delta/\| e^{-2f}\|_\infty}-Q_{\delta}
& = \prod_{\substack{[\ell]\in\pi(\hat{\mathcal L}_{\delta/\| e^{-2f}\|_\infty})\\ [\ell]\in\pi(\mathcal L_{\delta})}}\tr(\mathcal Hol^\nabla(\ell))\cdot\Big(\prod_{\substack{[\ell]\in\pi(\hat{\mathcal L}_{\delta/\| e^{-2f}\|_\infty})\\ [\ell]\notin\pi(\mathcal L_{\delta})}}\tr(\mathcal Hol^\nabla(\ell))-1\Big) \\
& \eqqcolon Q_{\delta} \cdot\left(X_\delta-1\right).
\end{align*}
We will show that $X_\delta$ converges toward $1$ in $L^2$; since $|Q_{ \delta}|\leq 1$, it will be sufficient to conclude the proof. To show that $X_\delta \to 1$ in $L^2$ amounts to show that $\mathbb{E}[X_\delta]\to 1$ and $\mathbb{E}[X_\delta^2]\to 1$, i.e. that 
$\log\mathbb{E}[X_\delta]\to 0$ and $\log\mathbb{E}[|X_\delta|^2]\to 0$, since the expectations here are real by symmetry, and positive by Campbell's formula.

On the one hand, using Campbell's formula, then Equation \eqref{eq:sym}, and finally the technical Proposition \ref{prop:tech}, we have 
\begin{align*}
0\leq -\log\mathbb{E}[X_\delta]&=\int_0^{  \delta } \int_M \frac{p_t(x,x)}{t} \mathbb{E}^m_{t,x,x}[ \mathbf1_{\hat T_W> \delta/\|e^{-2f}\|_\infty} (1-\tr (\mathcal Hol^\nabla(W)  ) ]\d x \d t \\
&\leq \int_0^{  \delta } \int_M \frac{p_t(x,x)}{t} \frac{1}{2} \mathbb{E}^m_{t,x,x}[ |1-\tr (\mathcal Hol^\nabla(W) ) |^2 ]\d x \d t \\
&\leq C\int_0^{ \delta }  \int_M\frac{t^{-2/2}}t t^2 \d x \d t= C\int_0^{  \delta }  \int_M \d x \d t   \underset{\delta\to 0}\longrightarrow 0.
\end{align*}
On the other hand, using $1-|z|^2\leq2\Re(1-z)$,
\begin{align*}
0\leq -\log\mathbb{E}[|X_\delta|^2]
&=\int_0^\infty\int_M\frac{p_t(x,x)}t\mathbb E^m_{t,x,x}\Big[\mathbf1_{T_W\leq\delta}\mathbf1_{\hat T_W>\delta/\|e^{-2f}\|_\infty}\big(1-\big|\tr(\mathcal Hol^\nabla(W))\big|^2\big)\Big]\mathrm dx\mathrm dt \\
&\leq \int_0^{\delta}\int_M\frac{p_t(x,x)}t\mathbb E^m_{t,x,x}\Big[2\Re\big(1-\tr(\mathcal Hol^\nabla(W))\big)\Big]\mathrm dx\mathrm dt \\
&\leq C\int_0^{\delta}\int_M t^0 \mathrm dx\mathrm dt \qquad\text{by Proposition \ref{prop:tech} and equation \eqref{eq:unifpbound}.}
\end{align*}
The integral is finite, hence converges to zero as $\delta\to0$, and Theorem \ref{th:conformal} is established.
\end{proof}

\appendix

\section{Continuity of holonomies over Brownian bridges}
\label{app:Econt}

We prove here continuity of the expectations of the holonomies over Brownian bridges, when we vary the endpoints and duration. Let us first recall basic properties of Brownian bridges: we work over $\widehat M$ closed, together with the usual base data of a metric $g$ on $\widehat M$, and the bundle data of a metric bundle $(E,h_E)$ endowed with a metric connection $\nabla$, all objects being smooth.

Given two vector bundles $F,F'$, respectively over $M,M'$, we write $F\boxtimes F'$ for their external tensor bundle, i.e. the vector bundle over $M\times M'$ whose fibre over $(x,y)$ is $F_x\otimes F'_y$. Such bundles appear naturally when considering holonomies, since the data of an endomorphism from $F_{\gamma_s}$ to $F_{\gamma_t}$, for $\gamma$ a curve with values in $M$, is equivalent to the data of an element of the fibre of $F\boxtimes F$ over $(\gamma_s,\gamma_t)$. To wit, set under $\mathbb P_{t,x,y}$ the parallel transport
\[ U^{TM}:s\mapsto\mathcal Hol^{\nabla^{TM}}_{0,st}(W) \]
of a frame over the Brownian bridge, seen as an element of the bundle $T^*M\boxtimes TM$ over $M\times M$, as well as the holonomy of $\nabla$ over the curve
\[ U^E:s\mapsto\mathcal Hol^\nabla_{0,st}(W), \]
which takes values in $E^*\boxtimes E$.
Denote by $\pi:T^*M\boxtimes TM\to M$ the projection to the second factor of the base space. It is known \cite[Theorem 5.4.4.]{Hsu} that $U^{TM}$ can be seen as the solution to an explicit SDE, while $U^E$ is defined as the solution to the implicit Stratonovich equation
\begin{align*} {\phantom{\mathrm (d^\nabla\text{ means }\nabla_{{}\circ\mathrm dW_{st}})}}\qquad\qquad\mathrm d^\nabla(U^E_sv)=0\qquad\qquad  \mathrm (d^\nabla\text{ means }\nabla_{{}\circ\mathrm dW_{st}}) \end{align*}
for every $v\in E_x$, made explicit in equation \eqref{eq:SDEHol}. Putting them together, we have for any reference connection $\nabla_0$%
\footnote{We use $\nabla_0$ to identify $\nabla$ with the section of a simple bundle, namely $(\nabla-\nabla_0)_x\in T^*_xM\otimes\operatorname{End}(E_x)$.}
\begin{equation}
\label{eq:SDEframe}
\left\{
\begin{aligned}
\mathrm dU^{TM}_s &= \sqrt tG_1(U^{TM}_s)\circ\mathrm du(w)_s + tG_2\big(U^{TM}_s,(\nabla p)_{t(1-s)}(W_s,y)\big)\mathrm ds, \\
\mathrm dU^E_s &= G_3\big((\nabla-\nabla_0)_{W_s},U^E_s\big){}\circ\mathrm dW_s
\end{aligned}
\right.
\end{equation}
for $G_1,G_2,G_3$ smooth independent of $(t,x,y)$ and $\nabla$, and $w$ a standard Brownian motion in $\mathbb R^d$  identified via an inessential isometry $u:\mathbb R^d\to T_xM$ to a Brownian motion $u(w)$ in $T_xM$.
In particular, thinking of $W$ as shorthand for $\pi(U^{TM})$, the dependence of $(U^{TM},U^E)$ with respect to $(t,x,y)$ can be abstracted as the dependence of the solution to some stochastic differential equation with respect to its coefficients and initial condition. Moreover, restricted to $s\in[0,1-\eta]$, the coefficients of these equations depend continuously on $(t,x,y)$ in the $\mathcal C^\infty$ topology; note that there is a singularity at $s=1$, because they involve the derivatives of the kernel $p$ for very small times.

Technically, the Stratonovich increment ${}\circ\mathrm du(w)_s$ depends on $x$ through the identification $u$. In our applications, we will take $x\mapsto u_x$ to be a smooth family of isometries, which must exist locally. This could be equivalently rephrased by taking $G_1$ to be a function of $(x,U^{TM}_s)$, although it is only well-defined locally in $x$. We leave this dependence implicit in our arguments, and claim that it does not impact the results.
\medskip

Let us go back to our original goal, using again the notations of section \ref{sec:not}.

\begin{lemma}
\label{le:Econt}
The function
\[ (t,x,y)\mapsto\mathbb E^m_{t,x,y}\big[\mathcal Hol^\nabla_{t,0}(W)\big] \]
is continuous over $(0,\infty)\times(\operatorname{int}M)^2$.
\end{lemma}

\begin{proof}
    Assume first that $M$ has no boundary.
    Fix $\alpha\in(1/3,1/2)$. We aim to find a collection of random processes $e^{t,x,y}:[0,1]\to E^*\boxtimes E$ indexed by $(t,x,y)$, defined on a single probability space, such that their distribution under the ambient probability measure $\mathbb P$ be that of
    \[ s\mapsto e^{-\int_0^{st} m(W_u)\d u}\mathcal Hol^\nabla_{0,st}(W) \]
    under $\mathbb P_{t,x,y}$, and such that the map
    \[ (t,x,y)\mapsto (e^{t,x,y}_s)_{s\in[0,1-\eta]}  \]
    be continuous as a map with values in $\mathcal C^\alpha([0,1-\eta],E^*\boxtimes E)$, for all $\eta>0$ small enough. Here the topology of the function space does not depend on the metric on $E\boxtimes E^*$ provided it is smooth, but we fix one arbitrarily and denote it, by abuse of notation, by $(\epsilon,\epsilon')\mapsto|\epsilon'-\epsilon|$. If these $e^{t,x,y}$ were to exist, then for every $M>0$, $s\in(1/2,1)$ and $\varepsilon>0$ satisfying the condition $2M(1-s)^\alpha\leq\varepsilon/3$, we would have by reflection invariance
    \begin{multline*}
    \mathbb P\big(\big|e^{t',x',y'}_1-e^{t,x,y}_1\big|\geq\varepsilon\big) \\
    \leq \mathbb P\big(\big|e^{t',x',y'}_{|[0,1/2]}\big|_{\mathcal C^\alpha}\geq M\big)
       + \mathbb P\big(\big|e^{t,x,y}_{|[0,1/2]}\big|_{\mathcal C^\alpha}\geq M\big)
       + \mathbb P\big(\big|e^{t',x',y'}_s-e^{t,x,y}_s\big|\geq\varepsilon/3\big),
    \end{multline*}
    so
    \[   \limsup_{(t',x',y')\to(t,x,y)}\mathbb P\big(\big|e^{t',x',y'}_1-e^{t,x,y}_1\big|\geq\varepsilon\big)
    \leq 2\mathbb P\big(\big|e^{t,x,y}_{|[0,1/2]}\big|_{\mathcal C^\alpha}\geq M\big) \]
    for all $M,\varepsilon>0$ by the continuity property above. Since the norm is finite almost surely we get the convergence of $e^{t',x',y'}_1$ toward $e^{t,x,y}_1$ in probability, hence by boundedness the continuity of
    \[ (t,x,y)\mapsto\mathbb E\big[e^{t,x,y}_1\big]=\mathbb E_{t,x,y}\big[e^{-\int_0^t m(W_s)\d s} \mathcal Hol^\nabla_{t,0}(W)\big]. \]
    
    Because $e$, or rather $(U^{TM},U^E,e)$, is solution to some extension of the Stratonovich equation \eqref{eq:SDEframe}, namely with the added
    \[ \mathrm de_s = e_s(U^E_s)^{-1}\circ\mathrm dU^E_s - t\cdot m(W_s)e_s\mathrm dt, \]
    the objective reduces to proving existence of a solution map that is continuous with respect to the Hölder topology for $e$.
    Although it must be possible to conclude using classical stochastic calculus, we were not able to come up with such a proof or find it in the literature, so we resorted to rough path machinery. Here is the result we use. Consider the Stratonovich stochastic differential equation
    \[ \mathrm dz_t = X_\theta(z_t)\circ\mathrm dw_s \]
    with initial value $p$ in some space manifold $N$. Here $t\mapsto w_t$ is a standard Brownian motion in $\mathbb R^d$, $t\mapsto z_t$ is a process with values in $N$, $\theta$ is an element of some manifold of parameters $\Theta$ and $X_\theta(z)$ is a linear map from $\mathbb R^d$ to $T_zN$. Then if $X$ is jointly smooth and the solutions are well-defined up to time 1 almost surely, there exists a coupling of all solutions $z=z^{p,\theta}$ such that the map $(p,\theta)\mapsto z^{p,\theta}$ is almost surely continuous from $N\times\Theta$ to $\mathcal C^\alpha([0,1],N)$. This is done by choosing $z$ the solution to the corresponding rough differential equation, with driving process a Brownian rough path---see e.g. \cite[Theorem 8.5; Theorem 9.1]{FrizHairer}.
    
    In our case, we get a coupling such that the map
    $(t,x,y)\mapsto e^{t,x,y}$
    is continuous over $(0,\infty)\times M\times M$, with values in $\mathcal C^\alpha([0,1-\eta],E^*\boxtimes E)$. Global existence is ensured by compactness of the spaces of isometries we consider. We stress again that rough path theory does not give us continuity for $\eta=0$, since the coefficients become ill-behaved; we do however have as much continuity as we needed.
\medskip
    
    This concludes the proof in the case $M$ has no boundary. For the general case, we assume as discussed in the introduction that $M$ is embedded in some closed manifold $\widehat M$ as a domain with a smooth boundary. Consider the processes $e^{t,x,y}$ constructed above as defined in $\widehat{M}$, and $W^{t,x,y}=\pi(e^{t,x,y})$. Then, what we need to show is that 
    \[ f:(t,x,y)\mapsto \mathbb{E}[ \mathbf{1}_{\Range(W^{t,x,y})\subset M}   e^{t,x,y}_1 ] \]
is continuous on $(0,\infty)\times(\operatorname{int}M)^2$.
Thus, it suffices to remark that, for $x,y,x',y'\in\operatorname{int}M$,
\begin{multline*}
|f(t,x,y)-f(t',x',y')| \\
\begin{aligned}
&\leq \mathbb{E}[|e^{t,x,y}_1-e^{t',x',y'}_1|]+ |\mathbb{E}[ e^{t',x',y'}_1( \mathbf{1}_{\Range(W^{t,x,y})\subset M}-\mathbf{1}_{\Range(W^{t',x',y'})\subset M}  ) ]|\\
&\leq \mathbb{E}[|e^{t,x,y}_1-e^{t',x',y'}_1|]+2 
|\mathbb{P}( \Range(W^{t,x,y})\subset M)-\mathbb{P}( \Range(W^{t',x',y'})\subset M)|\\
&\leq \mathbb{E}[|e^{t,x,y}_1-e^{t',x',y'}_1|]+2 
\left| \frac{ p^D_t(x,y)  }{p_t(x,y)}- \frac{ p^D_{t'}(x',y')  }{p_{t'}(x',y')}\right|,
\end{aligned}
\end{multline*}
where $p^D$ is the Dirichlet heat kernel associated with the Laplace--Beltrami operator on $M$.
Since the $e^{t,x,y}_1$ are uniformly bounded, they are uniformly integrable, and the argument above shows that the first summand must vanish as $(t',x',y')$ goes to $(t,x,y)$. By the continuity of the heat kernel and Dirichlet heat kernel, so does the second summand, and the proof is complete.
\end{proof}

\section{A multiplicative Campbell formula}
\label{app:Campbell}

We prove Theorem \ref{th:multiplicativeCampbell}. The following additive version is from \cite[Section 3.2]{Kingman}.

\begin{theorem}
\label{th:additiveCampbell}
Let $S$ be a measurable space such that the diagonal is measurable in $S\times S$.
Let $\mathcal{P}$ be a Poisson process on $S$ with intensity measure $\mu$, and let $f:S\to \mathbb{R}$ be measurable. The sum $\Sigma\coloneqq \sum_{X\in \mathcal{P}} f(X)$ is almost surely absolutely convergent if and only if
\[\int_S 1\wedge |f(x)| \mu(\d x) <\infty. \]
In this case, if
\[ \int_S |e^{f(x)}-1|\mu(\d x)<\infty, \]
then it holds that $\mathbb{E}[|e^\Sigma|]<+\infty$ and 
\[ 
\mathbb{E}[e^{ \Sigma}]= \exp\left(\int_S (e^{  f(x)}-1)\mu(\d x)\right).
\]
\end{theorem}

Let $S$, $\mathcal P$, $\mu$ and $g$ be as in Theorem \ref{th:multiplicativeCampbell} above.

Let us assume for now that $g$ never takes the value $-1$. Define over $\mathbb C\setminus\{-1\}$ the functions
\begin{align*} a:z&\mapsto\log|1+z|, & b:z&\mapsto\arg(1+z) \end{align*}
for $\arg:\mathbb C^*\to\left(-\pi,\pi\right]$ the complex argument. In some neighbourhood of $0$, $a(z)$ is of order $|z|$, and it is locally bounded. Accordingly, there exists a constant $C>0$ such that
\begin{align*}
      \int_S1\wedge\big|(a\circ g)(x)\big|\\mu(\mathrm dx)
&\leq \int_{\{|g|<1/2\}}\big|(a\circ g)(x)\big|\mu(\mathrm dx)
    + \int_{\{|g|\geq1/2\}}1\mu(\mathrm dx) \\
&\leq (3\vee C)\int_S|g(x)|\mu(\mathrm dx) < \infty, \end{align*}
and by the real additive Campbell formula, we know that the sum
\[ A\coloneqq\sum_{x\in\mathcal P}(a\circ g)(x) \]
is almost surely absolutely convergent. Following the exact same reasoning, the same holds for the sum
\begin{align*}
B&\coloneqq\sum_{x\in\mathcal P}(b\circ g)(x).
\end{align*}
Denoting by $\log_{\mathbb C}$ the complex logarithm from $\mathbb C^*$ to $\mathbb R+\mathsf i(-\pi,\pi]$, note that the sum
\[ \sum_{x\in\mathcal P}\log_{\mathbb C}(1+g(x)), \]
as the sum $A+\mathsf iB$, is absolutely convergent, so the product $\Pi$ is absolutely convergent.

Suppose that the integral
\[ \int_S\sup_{|w|<2}\left|\exp\Big(\big(a+wb\big)\big(g(x)\big)\Big)-1\right|\mu(\mathrm dx) \]
is finite, where the supremum ranges over complex numbers $w$ with $|w|<2$. Then we know that the function
\[ w\mapsto \int_S\left(\exp\Big(\big(a+wb\big)\big(g(x)\big)\Big)-1\right)\mu(\mathrm dx) \]
is well-defined and holomorphic over $B_0(2)$. Moreover, note that for all $\varepsilon>0$, $|w|\leq 2-\varepsilon$,
\[   \big|\exp(A+wB)\big|
\leq \exp\big(A+(2-\varepsilon)B\big) + \exp\big(A-(2-\varepsilon)B\big). \]
Since
\[ \int_S\left|\exp\Big(\big(a\pm(2-\varepsilon)b\big)\big(g(x)\big)\Big)-1\right|\mu(\mathrm dx)<\infty, \]
we know by the real additive Campbell formula that
\[ \mathbb E\left[\sup_{|w|\leq2-\varepsilon}\big|\exp\big(A+wB\big)\big|\right]<\infty, \]
and the function
\[ w\mapsto \mathbb E\left[\exp(A+wB)\right] \]
is well-defined and holomorphic over $B_0(2)$. Since we already know that the two functions above coincide for $w\in(-2,2)$, they must coincide throughout their domain, and in particular for $w=\mathsf i$ it yields
\begin{align*} \mathbb E[\Pi]
&= \mathbb E\big[\exp\big(A+\mathsf iB\big)\big] \\
&= \exp\left(\int_S\left(\exp\Big(\big(a+\mathsf ib\big)\big(g(x)\big)\Big)-1\right)\mu(\mathrm dx)\right)
 = \exp\left(\int_Sg(x)\mu(\mathrm dx)\right). \end{align*}

It remains to show the finiteness of the integral. The functions $a$ and $b$ are locally Lipschitz and vanish at zero, so there exists a constant $C>0$ such that for all $|z|<1/2$ and $|w|<2$,
\[   \left|\exp\Big(\big(a+wb\big)(z)\Big)-1\right| \leq C|z|. \]
If on the other hand $|z|\geq1/2$ with $z\neq-1$,
\begin{align*}
      \left|\exp\Big(\big(a+wb\big)(z)\Big)-1\right|
&\leq \exp\big(a(z)\big)\exp(|w|\pi)+1
    = |1+z|\exp(|w|\pi)+1 \\
&\leq (2|z|+|z|)\exp(|w|\pi) + 2|z| \\
&\leq 5\exp(|w|\pi)\cdot|z|.
\end{align*}
Up to taking a larger constant $C$, we find
\[   \sup_{|w|<2}\left|\exp\Big(\big(a+wb\big)\big(g(x)\big)\Big)-1\right|
\leq C|g(x)|, \]
and since the right-hand side is integrable, the left hand side must be as well. As discussed in the previous paragraph, this concludes the proof in the case where $g$ does not take the value $-1$.

If $g$ takes the value $-1$, set $\mathcal P_{=-1}$ and $\mathcal P_{\neq-1}$ the restrictions of $\mathcal P$ to the subsets $\{g=-1\}$ and $\{g\neq-1\}$ of $S$. These are independent Poisson processes, so in particular by the above the product
\[ \Pi_{\neq-1}\coloneqq\prod_{x\in \mathcal P_{\neq-1}}\big(1+g(x)\big) \]
is absolutely convergent almost surely, integrable, and
\[ \mathbb E[\Pi_{\neq-1}] = \exp\left(\int_Sg(x)\mu(\mathrm dx)\right). \]
Setting $N$ the cardinality of $\mathcal P_{=-1}$, we know that it is a Poisson random variable of (finite) parameter
\[ \mu(\{g=-1\}) = -\int_{\{g=-1\}}g(x)\mu(\mathrm dx), \]
and the product
\[ \Pi_{=-1}\coloneqq\prod_{x\in \mathcal P_{=-1}}\big(1+g(x)\big) \]
is finite. This means in particular that $\Pi=\Pi_{=-1}\cdot\Pi_{\neq-1}$ is almost surely absolutely convergent. Moreover, it is nothing but $0^N$, so the product is bounded, and by independence
\begin{multline*}
\mathbb E[\Pi]
= \mathbb E\big[\Pi_{=-1}\big]\cdot\mathbb E\big[\Pi_{\neq-1}\big]
= \mathbb P(N=0)\cdot\mathbb E\big[\Pi_{\neq-1}\big] \\
= \exp\left(\int_{\{g=-1\}}g(x)\mu(\mathrm dx)\right)\cdot\exp\left(\int_{\{g\neq-1\}}g(x)\mu(\mathrm dx)\right)
= \exp\left(\int_Sg(x)\mu(\mathrm dx)\right)
\end{multline*}
as announced.

\section*{Acknowledgments}
This version of the article has benefitted from many helpful comments and suggestions given by two anonymous reviewers, during the submission process to \emph{Annals of Probability.} We are very grateful, and are pleased to thank them for the numerous improvements they made possible.

\section*{Funding}
Pierre Perruchaud was partially supported by the Luxembourg National Research Fund -- grant O22/17372844/FraMStA.

Isao Sauzedde would like to thank the Isaac Newton Institute for Mathematical Sciences, Cambridge, for the support and hospitality during the programme \emph{Self-interacting processes} where work on this paper was undertaken. This work was supported by the EPSRC grants no EP/R014604/1 and EP/W006227/1.


\end{document}